\documentclass[final,3p,times]{elsarticle}
\usepackage{graphicx}
\usepackage{amssymb}
\usepackage{amsmath}
\usepackage{subfig}
\usepackage{graphicx}
\usepackage{url}
\usepackage{algorithm}
\usepackage{algpseudocode}
\usepackage{multirow}
\usepackage{wasysym}
\usepackage{marvosym}
\usepackage{color, xcolor, colortbl}
\usepackage{array, booktabs}
\usepackage{xspace}
\usepackage{bm}
\newcommand{\secref}[1]{Section \ref{#1}}
\newcommand{\fref}[1]{Fig.~\ref{#1}}
\newcommand{\tref}[1]{Table~\ref{#1}}

\renewcommand{\eqref}[1]{Eq.~(\ref{#1})}
\newcommand{\eref}[1]{(\ref{#1})}

\newcommand{\mat}[1]{\textrm{\textbf{#1}}}

\newcommand{\psif}{\psi (\bm{r}, \bm{\Omega}, E)}
\newcommand{\psifd}{\psi (\bm{r}, \bm{\Omega}', E')}

\definecolor{light}{rgb}{0.8,0.8,0.8}
\definecolor{medium}{rgb}{0.6,0.6,0.6}
\definecolor{dark}{rgb}{0.4,0.4,0.4}
\definecolor{darkmed}{rgb}{0.3,0.3,0.3}
\definecolor{darkest}{rgb}{0.2,0.2,0.2}
\definecolor{Black}{rgb}{0,0,0}
\definecolor{White}{rgb}{1,1,1}
\definecolor{lightpurple}{rgb}{0.78823,0.709803,0.74509}
\definecolor{lightpurpletext}{rgb}{0.788235,0.5529411,0.658823}
\definecolor{skyblue}{rgb}{0.80392,0.866666,0.92941}
\definecolor{skybluetext}{rgb}{0.61568627,0.7647058,0.913725}
\definecolor{darkgreen}{rgb}{0.3137254,0.458823,0.18431}
\definecolor{foliagegreen}{rgb}{0.188,0.415,0.105}
\definecolor{steelbluegrey}{rgb}{0.1961,0.2353,0.2392}
\definecolor{highlightblue}{rgb}{0.4078,0.6431,0.85}
\definecolor{matlabblue}{rgb}{0,0.2705,0.85}
\definecolor{darkred}{rgb}{0.8,0.1725,0}
\definecolor{fireenginered}{rgb}{0.505,0.1411,0}
\definecolor{darkpurple}{rgb}{0.6431,0.3137,0.8509}
\definecolor{gaylordpurple}{rgb}{0.416,0.204,0.549}
\definecolor{deludedorange}{rgb}{0.7409,0.4392,0}
\definecolor{darksalmon}{rgb}{0.9137,0.411,0.706}
\newcolumntype{a}{>{\columncolor{light}}c}
\newcolumntype{b}{>{\columncolor{skyblue}}c}
\begin{document}

\begin{frontmatter}


\title{A comparison of element agglomeration algorithms for unstructured geometric multigrid}
\author[AMCG]{S. Dargaville}
\ead{dargaville.steven@gmail.com}
\author[QM]{A.G. Buchan}
\author[AWE,AMCG]{R.P. Smedley-Stevenson}
\author[AMEC,AMCG]{P.N. Smith}
\author[AMCG]{C.C. Pain}
\address[AMCG]{Applied Modelling and Computation Group, Imperial College London, SW7 2AZ, UK}
\address[QM]{School of Engineering and Materials Sciences, Queen Mary University of London, E14 NS, UK}
\address[AWE]{AWE, Aldermaston, Reading, RG7 4PR, UK}
\address[AMEC]{ANSWERS Software Service, Jacobs, Kimmeridge House, Dorset Green Technology Park, Dorchester, DT2 8ZB, UK}
%~~~~~~~~~~~~~~~~~~~~~~~~~~~~~~~~~~~~
%~~~~~~~~~~~~~~~~~~~~~~~~~~~~~~~~~~~~
\begin{abstract}
This paper compares the performance of seven different element agglomeration algorithms on unstructured triangular/tetrahedral meshes when used as part of a geometric multigrid. Five of these algorithms come from the literature on AMGe multigrid and mesh partitioning methods. The resulting multigrid schemes are tested matrix-free on two problems in 2D and 3D taken from radiation transport applications; one of which is in the diffusion limit. In two dimensions all coarsening algorithms result in multigrid methods which perform similarly, but in three dimensions aggressive element agglomeration performed by METIS produces the shortest runtimes and multigrid setup times. 
\end{abstract}
%~~~~~~~~~~~~~~~~~~~~~~~~~~~~~~~~~~~~

%~~~~~~~~~~~~~~~~~~~~~~~~~~~~~~~~~~~~
%~~~~~~~~~~~~~~~~~~~~~~~~~~~~~~~~~~~~

\end{frontmatter}
%~~~~~~~~~~~~~~~~~~~~~~~~~~~~~~~~~~~~
%~~~~~~~~~~~~~~~~~~~~~~~~~~~~~~~~~~~~
%~~~~~~~~~~~~~~~~~~~~~~~~~~~~~~~~~~~~
%-----------------------------
\section{Introduction}
\label{sec:Introduction}
All multigrid methods rely on some form of coarsening and for many mesh-based discretisations of PDEs, element agglomeration can be used to generate coarse grids. The use of element coarsenings is common in geometric multigrids and finite volume multigrid schemes \cite{dervieux_multigrid_1988, Mavriplis1996, Thomas2011, Nishikawa}, but its use has been growing in algaebraic methods (AMGe-type methods), along with upscaling techniques. This is based on the idea that access to the underlying finite element matrices, in addition to the assembled matrix can help improve the methods. Please see \cite{Venkatakrishnan1993, Henson2001, Brezina2001,  jones_amge_2001, Vassilevski2002, Chartier2003, wabro_amgecoarsening_2006, lashuk_versions_2008, Brezina2012, Kalchev2016, Botti2017} for examples of element agglomeration in these applications.  

We use element agglomeration algorithms to construct multigrid methods to solve problems involving the Boltzmann transport, in applications areas such as radiation transport, spectral wave and kinetic modelling. Boltzmann transport problems are difficult to solve in part because of their very large phase-space; radiation transport for example have a seven-dimensional phase space, with three spatial dimensions, two angular, one energy and time. The iterative methods used to solve these problems therefore must be entirely matrix-free, where even the unassembled element matrices are too large to form; often there is not enough memory to even store a single copy of the solution. This precludes many matrix-based multigrid approaches like AMG, AMGe or $\rho$-AMGe. 

As such, we use element agglomeration as part of a geometric multigrid on our spatial mesh applied matrix-free. This is based on the work of \cite{buchan_sub-grid_2012}, who use both a CG and DG spatial mesh as part of a sub-grid scale finite element method. Element agglomeration is then used to define geometrically how the CG and DG spatial tables are Galerkin projected onto lower grids in a agglomerate local fashion, to form coarse spatial tables. This is very similar to AMGe methods, but without the sophisticated operators derived from the unassembled element matrices. As mentioned above, we do not have the memory to construct the local element matrices and in particular, given the number of DOFs at each spatial node, inverting/decomposing the element matrix would be very expensive, even if we could store the resulting operators. 

We therefore construct simple geometric operators which prolong/restrict the DOFs on a spatial node in the same way (i.e., block operators with size that depends only on the spatial dimensions). This allows us to apply both the operators and matrix-vector products on all levels matrix-free, ensuring low memory use. We should note that if a structured spatial grid is used and the element agglomeration results in structured lower grids, this is equivalent to a simple rediscretisation on the coarser grid. 

Element agglomeration can be trivial on structured grids, but generating ``good'' agglomeration algorithms on unstructured grids can be difficult. Given the algebraic nature of AMGe algorithms, this information can also be used to affect the coarsening of the spatial mesh itself, e.g., \cite{jones_amge_2001} who use estimates of interpolation accuracy to inform the coarsening. Again given the lack of element matrices, we focus on a purely geometric coarsening in this work. The existing literature described above uses a range of different algorithms to generate coarsenings and it can be difficult to determine which would be optimal. 

The work of \cite{wabro_amgecoarsening_2006} compared the performance of an AMGe solver (in which they constructed energy-minimising operators) which used the Jones \cite{jones_amge_2001}, RGB, METIS \cite{karypis_metis-unstructured_1995} (in a bottom-up fashion) and MGridGen \cite{Moulitsas2001} algorithms to construct agglomerates. The authors solved both a convection-diffusion problem and the Onseen equation and compared the solve times along with several common metrics (e.g., operator and grid complexities). They showed considerable difference between the performance of the methods, with MGridGen and METIS both performing well. 

There is however one key parameter which determines the behaviour of many coarsening algorithms, and that is the desired number of fine elements in an agglomerate. \cite{wabro_amgecoarsening_2006} however only examined the performance of METIS in 2D when given two different ``target'' agglomerate sizes, namely 4 and 5. Not all coarsening algorithms require a target size (e.g., \citet{jones_amge_2001}), but for those that do it is crucial to the performance of the method. The target size on a structured cartesian grid is trivial, namely 2$^d$, where $d$ is the number of spatial dimensions. When extending to unstructured grids however, the desired size is harder to determine given the increased connectivity of unstructured meshes and the desire to keep the operator complexity low. 

We must note that this work is not focused on building an optimal multigrid for Boltzmann transport problems. We use very simple operators in this work, given our matrix-free constraints. Radiation transport problems range from very diffuse problems (which in the limit reduce to the diffusion equation, which is very suited to traditional multigrid methods), to multi-material problems with hyperbolic regions, which are very difficult to solve without sophisticated operators. 

Hence in this work, we compare the performance of seven different coarsening algorithms in forming coarse grids for our geometric, matrix-free multigrid, applied to the Boltzmann transport equation on unstructured triangular/tetrahedral meshes. In order to help extend the generality of this work outside of the Boltzmann transport community, we deliberately include a test problem within the diffusive limit, which helps give insight into the performance of these coarsening algorithms on more traditional elliptic operators. We also include a more challenging problem away from the diffusive limit, in both two and three spatial dimensions. Furthermore, we examine the choice of key parameters like the target agglomerate size and investigate the role of ``aggressive coarsenings'' from AMG methods in the context of element agglomeration on unstructured grids. This paper therefore aims to provide a practical guide to the choice of element coarsening algorithms for unstructured, geometric multigrids along with their resulting parameters. 
%-----------------------------
\section{Boltzmann transport equation}
\label{sec:Boltzmann}
In this work we solve the steady-state Boltzmann Transport Equation (BTE), which in the form given below governs the behaviour of neutral particles (e.g., neutrons) moving through a medium in radiation transport applications. These particles are advected outward from sources in all directions and can be absorbed or change energy/direction by the material they are moving through. We can therefore write the first-order integro-differential form of the BTE at steady state as
\begin{equation}
\bm{\Omega} \cdot \nabla \psif + \Sigma_\textrm{t} \psif = \int_{\bm{\Omega}'} \int_{E'} \Sigma_\textrm{s} (\bm{r}, \bm{\Omega}' \rightarrow \bm{\Omega}, E' \rightarrow E) \psifd \textrm{d}E' \textrm{d}\bm{\Omega}' + S_\textrm{e}(\bm{r}, \bm{\Omega}, E),
\label{eq:bte}
\end{equation}
where $\psif$ is the angular flux in direction $\bm{\Omega}$, at spatial position $\bm{r}$ and at energy $E$. The parameters $\Sigma_\textrm{t}$ and $\Sigma_\textrm{s}$ describe the material properties, giving the total and scatter macroscopic cross sections, respectively. Finally we include a source term $S_\textrm{e}$, which in this work can vary in space, but is constant in angle to model a material releasing particles isotropically in angle. We also assume only isotropic scattering.

Equation \ref{eq:bte} has a 6-dimensional phase-space involving three spatial, two angular and one energy dimension. We should note the two important limits of \eref{eq:bte}. The first is if $\Sigma_\textrm{s} \rightarrow \infty$ then \eref{eq:bte} reduces to a diffusion equation in space, as particles scatter heavily (i.e., there is no angular variation) and the mean free path is very small. Conversely, if $\Sigma_\textrm{s},\Sigma_\textrm{t} \rightarrow 0$, then \eref{eq:bte} describes particles which propogate an infinite distance without interacting (i.e., the mean free path is infinite). We simplify \eref{eq:bte} here by assuming only one energy for all the particles, hence reducing the phase-space to five dimensions. 

We must use a stabilised discretisation scheme in space, given the transport operator in \eref{eq:bte}. We use the sub-grid scale (SGS) finite element scheme described by \cite{buchan_inner-element_2010}; please see this reference for more details. To discretise in angle, we use the standard P$_n$ approach \cite{Lewis1993}, which is a spectral method on the sphere using spherical harmonic basis functions; we restrict the order in this work to P$_1$, giving 3 angular DOFs per spatial node in 2D, and 4 in 3D. 
%-----------------------------
\section{Multigrid method}
\label{sec:Multigrid method}
Before discussing the coarsening algorithms we examine in this work, we first describe our iterative method; we use our multigrid as a preconditioner for FGMRES(30) \cite{saad_flexible_1993}. Assuming a sequence of coarse grids have been formed by element agglomeration, we must chose our smoothers and interpolation/restriction operators. Traditionally in radiation transport a common smoother choice is a Gauss-Seidel (known as a ``sweep''), as the use of structured grids with DG-FEM in space and a discrete ordinates discretisation (known as S$_n$) in angle results in the streaming/removal part of the discretised transport equation (i.e., the LHS in \eref{eq:bte})  having lower triangular structure. This pure-advection problem can therefore be solved exactly in one Gauss-Seidel iteration, making a powerful smoother for \eref{eq:bte} (e.g., see \cite{Manteuffel1994, Oliveira1998b, Lee2010a, Kanschat2014, Densmore2016}).

On unstructured grids however, it can be difficult to perform optimal sweeps, particularly in parallel and that motivates our continuing work on investigating alternate smoothers. In this work we use a Jacobi-preconditioned GMRES(3) method as our smoother (as do \cite{buchan_sub-grid_2012, Dargaville2019, Dargaville2019a}), with three iterations on each level as the pre and post smoother. These strong smoothers help make up for our simple operators. 

Recently however, interpolation/restriction operators have been developed for the BTE that when used as part of an AMG method for solving the discretised streaming/removal terms lead to highly effective multigrid algorithms \cite{Manteuffel2019, Manteuffel2018}, even without Gauss-Seidel smoothers. Currently the use of these operators relies on a specific iterative scheme (``source iteration'' aka a Richardson iteration) which solves the streaming/removal problem separately to the scattering, rather than using a multigrid method on the entirety of \eref{eq:bte}. The ability to use this source iteration is precluded in many multi-physics Boltzmann applications as it relies on the mapping between the underlying S$_n$ quadrature and Legendre moments. 

Given this, we apply our multigrid method to all the terms in \eref{eq:bte} and hence use very simple interpolation/prolongation operators (largely following \cite{buchan_sub-grid_2012}). Our prolong/restrict is applied identically across each angular DOF. If we consider coarse and fine spatial nodes, for coarse nodes that exist on the finer grid, we use injection to prolong their values; for fine nodes on a coarse face, they are prolonged using an average of the coarse nodes on the coarse face; finally for fine nodes on the interior of coarse elements, they are formed by averaging the coarse nodes they are directly connected to on the fine element. 

If there are no coarse nodes on the fine element, we progressively expand our search out to neighbouring fine elements connected (in the finite element sense) to the original fine element. If we have very small agglomerates this search is rarely triggered, but for large agglomerates this helps reduce our operator complexity by keeping prolongation as local as possible. If we are forced to expand the search radius more than 20 times (on a structured quad grid in 2D with an optimal coarsening, expanding the neighbour 5 times would mean an agglomerate larger than 121 fine elements), we then simply take the average of all coarse nodes on the agglomerate in order to stop the search from becoming prohibitively expensive. This is only triggered on the largest of agglomerates, that are not shaped like circles/spheres (e.g., long, thin agglomerates). We set our restrictor to be the transpose of our prolongator. 

Our coarse grid finite element spatial tables are formed from element-local projections of the fine tables, and the material properties on coarse elements are formed from volume weighted averages of those on the finer grids. 

Given these constraints, the performance of the element-based coarsening algorithm on an unstructured mesh becomes crucial. It seems natural to try and emulate the coarsening of a structured grid (indeed many of the existing algorithms return perfect coarsenings on structured grids), with well-shaped (no ``holes'') compact agglomerates (low surface area to volume ratio) with low element connectivity (which corresponds to a small edge-cut when partitioning the dual-graph) and containing a fixed number of fine elements. Unfortunately, balancing these points on an unstructured grid while retaining acceptable grid-complexities can be challenging. 
%~~~~~~~~~~~
\section{Coarsening algorithms}
\label{sec:Coarsening algorithms}
%~~~~~~~~~~~
In general, any mesh partitioning tool can be used to generate an element-based coarsening. There are several common algorithms used in the literature, which are briefly described here. We have tried to keep the description of each coarsening algorithm simple (for more detail or more precise notation please see the original papers). On a set multigrid level, we use the notation $\mathcal{A}$ for agglomerates of elements (or coarse elements), $\mathcal{E}$ for individual elements (or fine elements), $f$ and $g$ for faces on elements, $e$ and $d$ for edges on faces (only in 3D), $w$ for weights applied to each face/edge and $\mathcal{N}$ for neighbours of an element (two elements are neighbours if they share a face).

The \citet{jones_amge_2001} algorithm (see Algorithm \ref{alg:jones}) constructs agglomerates by assigning weights to faces and constructs well-shaped agglomerates in 2D. \citet{kraus_agglomeration-based_2004} modified the algorithm of \citet{jones_amge_2001} by also weighting edges in order to improve the shape of agglomerates in 3D (see Algorithm \ref{alg:kraus}; we should note the Kraus algorithm also modifies the selection of a maximally-weighted face so the two algorithms will not produce the same coarsening in 2D). The RGB algorithm of \citet{wabro_amgecoarsening_2006} (see Algorithm \ref{alg:rgb}) coarsens aggressively based on local neighbours, producing well shaped agglomerates. The library MGridGen \cite{moulitsas_multilevel_2001} (see Algorithm \ref{alg:mgridgen}) solves an optimisation problem for agglomerates based on minimising metrics that include the surface area to volume ratio of agglomerates across the mesh (we use the default options). The graph-partitioning library METIS/PARMETIS \cite{karypis_metis-unstructured_1995} can also be used (in either a top-down or bottom-up fashion \cite{wabro_amgecoarsening_2006}) to produce an element-based coarsening (see Algorithm \ref{alg:metis}; we use it top-down, from the finest to coarsest mesh), minimising the edge-cuts of the dual-graph. We also set the option \textsc{METIS\_OPTION\_CONTIG} to try and force contiguous partitions.

The algorithms listed above vary in sophistication, but we can also construct very simple algorithms. Algorithm \ref{alg:node}, which we denote as the ``node'' algorithm, simply picks fine nodes on the mesh and agglomerates all elements connected to it (this is very similar to early finite volume agglomeration schemes). We also consider a ``greedy'' algorithm (described in Algorithm \ref{alg:greedy}), that selects fine elements and adds neighbouring elements (with no preference given to ordering) until reaching a desired agglomerate size.

In total, this gives 7 different coarsening algorithms. To give a sense of the agglomerates produced by each algorithm, Tables \ref{tab:coarse_unstruct} and \ref{tab:coarse_unstruct3D} show the coarsenings produced on a triangular meshing of the unit square in 2D and a tetrahedral meshing of the unit cube in 3D, respectively. In 2D, most algorithms produce good looking agglomerates (with the possible exception of the Greedy algorithm, which does not attempt to optimise the shape of agglomerates in any way). In 3D, it is more difficult to judge the shape of each agglomerate, though it should be noted that many of the algorithms produce agglomerates with ``holes''. We should note that several of the algorithms described produce optimal or near-optimal coarsenings on structured grids, and this can be a very important property to preserve. We do not explicitly examine the coarsenings on structured grids, but note that all but the simplest algorithms (e.g., greedy) in this work can produce very good structured coarsenings. 

Most of these algorithms can also be easily extended, for example to preserve material boundaries across coarse agglomerates as much as possible. Indeed this can be considered one of the benefits of element coarsening, in that this information is easily accessible. We did test using the Jones, Greedy, METIS and MGridGen algorithms with material boundaries preserved on coarse grids, but we did not see a large difference in our runtimes. The relative difference in material properties however in this work are not large; we leave examining this to future work and note the wealth of multigrid literature on the benefits of preserving material properties.
%~~~~~~~
\begin{algorithm}[ht]
\footnotesize
\caption{\citet{jones_amge_2001} coarsening algorithm}\label{alg:jones}
\begin{algorithmic}[1]
\State Set $w(f)=0$ for all faces $f$
\State \textbf{Global search}: Find a face $f$ with maximal $w(f) \geq 0$. Create a new agglomerate $\mathcal{A}$ \label{global_search}
\State Set $\mathcal{A} = \mathcal{A} \cup \mathcal{E}_1 \cup \mathcal{E}_2$, where $\mathcal{E}_1 \cap \mathcal{E}_2 = f$ and set $w_\textrm{max} = w(f)$, $w(f) = -1$ \label{step_1}
\State Increment $w(f_1)$ for all faces $f_1$ that are neighbours of $f$ and $w(f_1) \neq -1$
\State Increment $w(f_2)$ for all faces $f_2$ that are neighbours of $f$, with $f_2$ and $f$ sharing an element and $w(f_1) \neq -1$
\State Choose the face $g$ with maximal weight from the neighbours of $f$
\If{$w(g) \geq w_\textrm{max}$}
\State Set $f=g$ and goto \ref{step_1}
\Else
\State $\mathcal{A}$ is complete. Set the weight of every face in $\mathcal{A}$ to -1 and goto \ref{global_search}
\EndIf
\State Cleanup the coarsening
\end{algorithmic}
\end{algorithm}
%~~~~~~~

%~~~~~~~
\begin{algorithm}[ht]
\footnotesize
\caption{\citet{kraus_agglomeration-based_2004} coarsening algorithm}\label{alg:kraus}
\begin{algorithmic}[1]
\State Set $w(f)=0$ for all faces $f$, $w(e)=0$ for all edges $e$ and case=``edge''
\State \textbf{Global search}: Find an edge $e$ with maximal $w(e) \geq 0$. If there are no edges remaining ($w(e)<0$), set case=``face'' and find a face $f$ with maximal $w(f) \geq 0$. Create a new agglomerate $\mathcal{A}$ \label{global_search_kraus}
%~~
%Step1
\If{case $==$ ``edge''} \label{step_1_kraus}
\State Set $\mathcal{A} = \mathcal{A} \cup \mathcal{E}_e$, where $\mathcal{E}_e$ are the elements that share $e$ and set $w_\textrm{$e$,max} = w(e)$, $w(e) = -1$ 
\Else
\State Set $\mathcal{A} = \mathcal{A} \cup \mathcal{E}_1 \cup \mathcal{E}_2$, where $\mathcal{E}_1 \cap \mathcal{E}_2 = f$ and set $w_\textrm{$f$,max} = w(f)$, $w(f) = -1$
\EndIf
%~~
%Step2
\If{case $==$ ``edge''}
\State Increment $w(e_1)$ for all edges $e_1$ that are neighbours of $e$ and $w(e_1) \neq -1$
\Else
\State Increment $w(f_1)$ for all faces $f_1$ that are neighbours of $f$ and $w(f_1) \neq -1$
\EndIf
%~~
%Step3
\If{case $==$ ``edge''}
\State Increment $w(e_2)$ for all edges $e_2$ that are neighbours of $e$, with $e_2$ and $e$ sharing a face and $w(e_2) \neq -1$
\Else
\State Increment $w(f_2)$ for all faces $f_2$ that are neighbours of $f$, with $f_2$ and $f$ sharing an element and $w(f_2) \neq -1$
\EndIf
%~~
%Step4 - edge only
\If{case $==$ ``edge''}
\State Increment $w(f_3)$ for all faces $f_3$ that are neighbours of any face that has edge $e$ on it,  $w(f_3) \neq -1$
\EndIf
%~~
%Step5
\If{case $==$ ``edge''}
\State Choose the edge $d$ with maximal weight from the neighbours of $e$, with $d$ and $e$ sharing a face
\If{$w(d) \geq w_\textrm{$e$,max}$}
\State Set $e=d$ and goto \ref{step_1_kraus}
\Else
\State $\mathcal{A}$ is complete. Set the weight of every face and edge in $\mathcal{A}$ to -1 and goto \ref{global_search_kraus}
\EndIf
\Else
\State Choose the face $g$ with maximal weight from the neighbours of $f$, with $g$ and $f$ sharing an element
\If{$w(g) \geq w_\textrm{$f$,max}$}
\State Set $f=g$ and goto \ref{step_1_kraus}
\Else
\State $\mathcal{A}$ is complete. Set the weight of every face $\mathcal{A}$ to -1 and goto \ref{global_search_kraus}
\EndIf
\EndIf
\State Cleanup the coarsening
\end{algorithmic}
\end{algorithm}
%~~~~~~~

%~~~~~~~
\begin{algorithm}[ht]
\footnotesize
\caption{RGB coarsening algorithm of \citet{wabro_amgecoarsening_2006}}\label{alg:rgb}
\begin{algorithmic}[1]
\While{there are uncoloured elements}
\State Select an uncoloured element $\mathcal{E}$ randomly and colour it black
\State Colour all uncoloured or grey neighbours of $\mathcal{E}$ red
\State Colour all uncoloured elements which neighbour red elements grey
\EndWhile
\State Black elements plus surrounding red elements form agglomerates
\State Grey elements are appended to the agglomerates where they ``fit best''
\State Cleanup the coarsening
\end{algorithmic}
\end{algorithm}
%~~~~~~~

%~~~~~~~
\begin{algorithm}[ht]
\footnotesize
\caption{Node-based coarsening algorithm}\label{alg:node}
\begin{algorithmic}[1]
\For{all the unused nodes not on a boundary}
\State Randomly select a node
\State Create a new agglomerate $\mathcal{A}$
\State Get the elements that share the node and add them to $\mathcal{A}$
\State Get all the nodes and elements in $\mathcal{A}$ and call them used 
\EndFor
\State Now repeat for all the unused nodes on a boundary
\State Cleanup the coarsening
\end{algorithmic}
\end{algorithm}
%~~~~~~~

%~~~~~~~
\begin{algorithm}[ht]
\footnotesize
\caption{Greedy coarsening algorithm}\label{alg:greedy}
\begin{algorithmic}[1]
\State Set desired agglomerate size, $s$
\While{there are unused elements}
\State Randomly select an unused element $\mathcal{E}$ and add it to a new agglomerate $\mathcal{A}$ (call $\mathcal{E}$ used)
\State Insert the unused neighbours of $\mathcal{E}$  into a set $\mathcal{N}$
\While{the length of $\mathcal{N}$ != 0} \label{neigh_while}
\State Set $\mathcal{E}_n$ to the first element in $\mathcal{N}$
\State Remove $\mathcal{E}_n$ from $\mathcal{N}$ and add it to $\mathcal{A}$ (call $\mathcal{E}_n$ used)
\If{length $\mathcal{N} ==s$}
\State Flush $\mathcal{N}$ and exit \ref{neigh_while} 
\EndIf
\State Insert the unused neighbours of $\mathcal{E}_n$ into $\mathcal{N}$
\EndWhile
\EndWhile
\State Cleanup the coarsening
\end{algorithmic}
\end{algorithm}
%~~~~~~~

%~~~~~~~
\begin{algorithm}[ht]
\footnotesize
\caption{METIS 5.1 \cite{karypis_metis-unstructured_1995} coarsening algorithm}\label{alg:metis}
\begin{algorithmic}[1]
\State Construct the dual-graph of the mesh
\State Build integer area/volume weights from the dual-graph by scaling relative face areas and volumes
\State Set desired agglomerate size, $s$
\State Set the number of partitions as: int(number of elements / $s$)
\If{number of partitions $>$ 8}
\State Call METIS\_PartGraphKway
\Else
\State Call METIS\_PartGraphRecursive
\EndIf
\State Cleanup the coarsening
\end{algorithmic}
\end{algorithm}
%~~~~~~~

%~~~~~~~
\begin{algorithm}[ht]
\footnotesize
\caption{MGridGen 1.0 \cite{moulitsas_multilevel_2001} coarsening algorithm}\label{alg:mgridgen}
\begin{algorithmic}[1]
\State Construct the dual-graph of the mesh
\State Set desired agglomerate size, $s$
\State Set the minimum and maximum agglomerate size to $s$
\State Call the MGridGen library
\State Cleanup the coarsening
\end{algorithmic}
\end{algorithm}
%~~~~~~~

%~~~~~~~~~~~
\subsection{Cleanup}
\label{sec:Cleanup}
All element-based coarsening algorithms require some kind of cleanup stage to ensure that the resulting agglomerates represent sensible coarse elements. Ideally, a coarsening algorithm would require no cleanup, though in practice we have found an effective cleanup stage to be vital in producing a robust multigrid method (particularly in 3D). Some common circumstances that require cleanup can include:
\begin{enumerate}
\item Unused fine elements. We append unused elements to existing agglomerates that share the most faces. If there are multiple agglomerates that share the same number of faces with the unused fine element, we chose the agglomerate with the smallest number of fine elements.
\item Unused fine elements surrounded by only unused fine elements. We identify the unused fine elements that neighbour existing agglomerates, append those to the agglomerate, then repeat for any unused fine elements that remain.
\item Agglomerates made up of unconnected fine elements or groups of fine elements (i.e., fine elements that are only connected by nodes, edges or are completely non-contiguous). If this occurs, we find the smallest group(s) of unconnected fine elements and append them to other agglomerates that they share the most faces with. 
\item Completely enclosed agglomerates (i.e., one agglomerate completely covers another). This is more common in 3D (particularly in corners of the domain) and can result in a lack of coarse nodes on an agglomerate. To prevent this, we simply combine the two agglomerates.
\item Faces/edges that lose sensible interpretations on lower grids. For example, the \citet{kraus_agglomeration-based_2004} algorithm uses both faces and edges in 3D to perform a coarsening. Computing edges on lower grids can sometimes be problematic, resulting in ``edge loops'' (edges that connect to themselves) or ``hanging edges'' (edges that are not continuously connected), which must be dealt with, often by simply excluding them from the list of coarse edges, or breaking them into multiple coarse edges.
\end{enumerate}

As such, all coarsening algorithms shown in \secref{sec:Coarsening algorithms} include a cleanup stage once the agglomeration is completed on each level.
%~~~~~~~~~~~
%~~~~~~~~~~~
%~~~~~~~~~~~
\subsection{Coarse face/edge/node selection}
\label{sec:Coarse node selection}
Now that we have a number of agglomeration algorithms, we need to determine the coarse faces, edges and nodes on the coarse elements. In order to select coarse faces, we follow \cite{jones_amge_2001, kraus_agglomeration-based_2004} and define the coarse faces as the intersection of fine faces between different agglomerates (this results in two faces between two neighbouring agglomerates), or between an agglomerate and the boundary of the domain. Depending on the coarsening algorithm used, we can also compute coarse edges and finally coarse nodes. If using any algorithm but \citet{kraus_agglomeration-based_2004}, we skip computing edges and determine the coarse nodes directly, by counting the number of coarse faces each fine node is in, and comparing it to the number of coarse elements that node is contained in. Given the number of spatial dimensions $d$, a fine node is considered coarse if the number of coarse faces $> 2^{d-2}$ number of coarse elements (i.e., coarse nodes are on the vertices of coarse elements). 

When using the \citet{kraus_agglomeration-based_2004} algorithm, we first compute coarse edges by using the intersection of fine edges between coarse faces. We should note that this can be very expensive (as edges are connected to many faces). Coarse node selection is then given by the intersection of fine nodes on a given coarse edge. For more details, please see \citet{kraus_agglomeration-based_2004} who give a precise definition of the topological relations used to create these coarse topologies.
%~~~~~~~~~~~
%~~~~~~~~~~~
%\begin{table}[htp]
%\begin{center}
%\begin{tabular}{c c c c c}
%\toprule
%Alg. & Level 0 & 1 & 2 & 3 \\
%\midrule
%1, 2, 3 & \includegraphics[height=2.5cm]{figures/square_16_level_0} &\includegraphics[height=2.5cm]{figures/square_16_m_1_level_1}&\includegraphics[height=2.5cm]{figures/square_16_m_1_level_2}&\includegraphics[height=2.5cm]{figures/square_16_m_1_level_3}\\
%4 & \includegraphics[height=2.5cm]{figures/square_16_level_0} &\includegraphics[height=2.5cm]{figures/square_16_m_4_level_1}&\includegraphics[height=2.5cm]{figures/square_16_m_4_level_2}&\\
%5 & \includegraphics[height=2.5cm]{figures/square_16_level_0} &\includegraphics[height=2.5cm]{figures/square_16_m_5_level_1}&\includegraphics[height=2.5cm]{figures/square_16_m_5_level_2}&\includegraphics[height=2.5cm]{figures/square_16_m_5_level_3}\\
%6 & \includegraphics[height=2.5cm]{figures/square_16_level_0} &\includegraphics[height=2.5cm]{figures/square_16_m_6_level_1}&\includegraphics[height=2.5cm]{figures/square_16_m_6_level_2}&\includegraphics[height=2.5cm]{figures/square_16_m_6_level_3}\\
%\bottomrule
%\end{tabular}
%\end{center}
%\label{fig:coarse_struct}
%\caption{Coarsening algorithms applied to a small (256 element) structured grid in 2D}
%\end{table}

\begin{table}[p]
\begin{center}
\begin{tabular}{>{\centering\arraybackslash} m{1.75cm} >{\centering\arraybackslash} m{2.5cm} >{\centering\arraybackslash} m{2.5cm} >{\centering\arraybackslash} m{2.5cm} >{\centering\arraybackslash} m{2.5cm} >{\centering\arraybackslash} m{2.5cm} >{\centering\arraybackslash} m{2.5cm}}
\toprule
Alg. & Orig. mesh & Level 1 & 2 & 3 & 4 & 5 \\
\midrule
Jones \cite{jones_amge_2001} & \includegraphics[height=2.5cm]{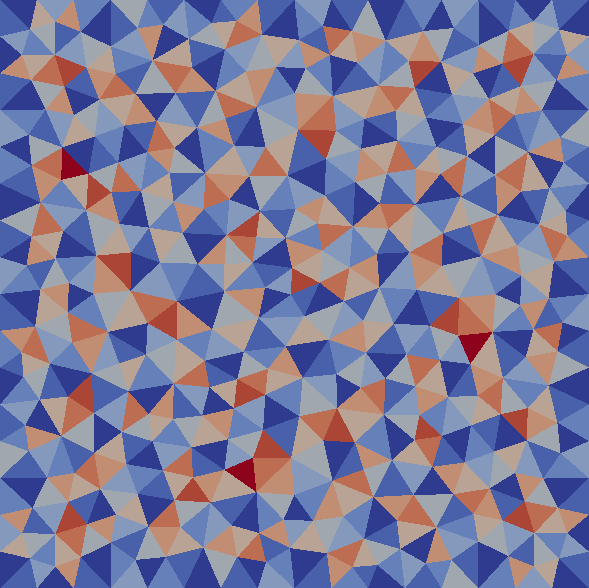} &\includegraphics[height=2.5cm]{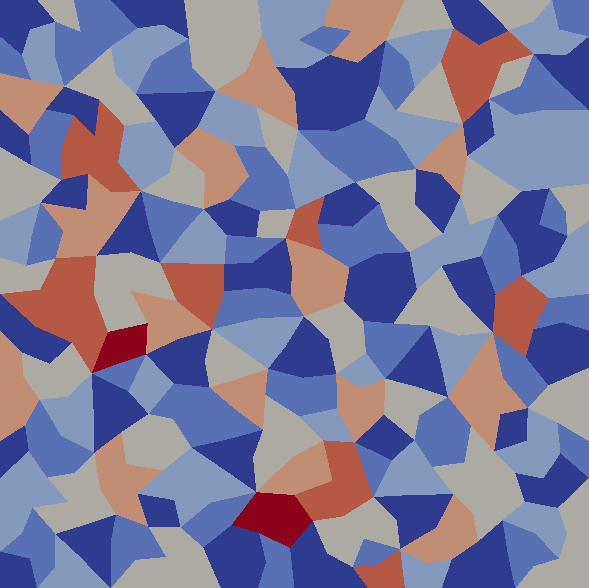}&\includegraphics[height=2.5cm]{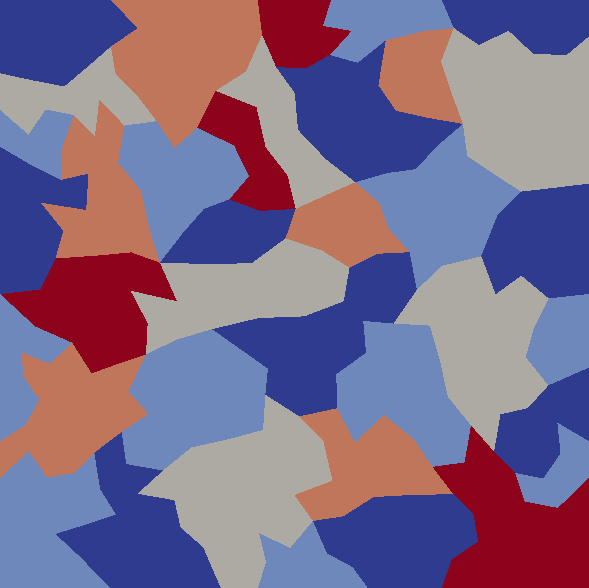}&\includegraphics[height=2.5cm]{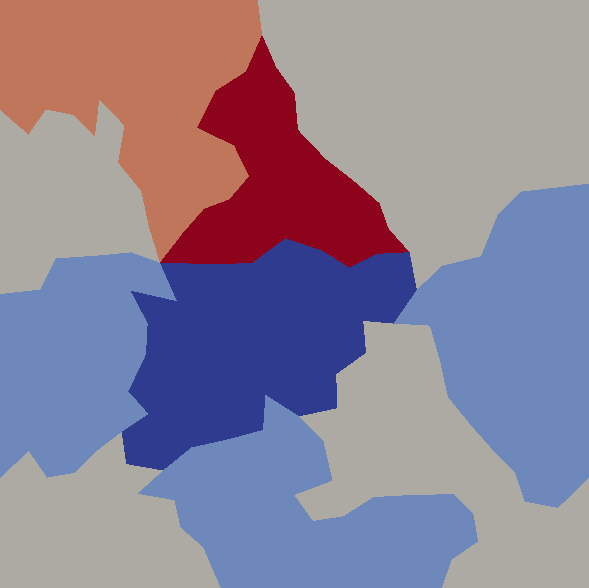}&\includegraphics[height=2.5cm]{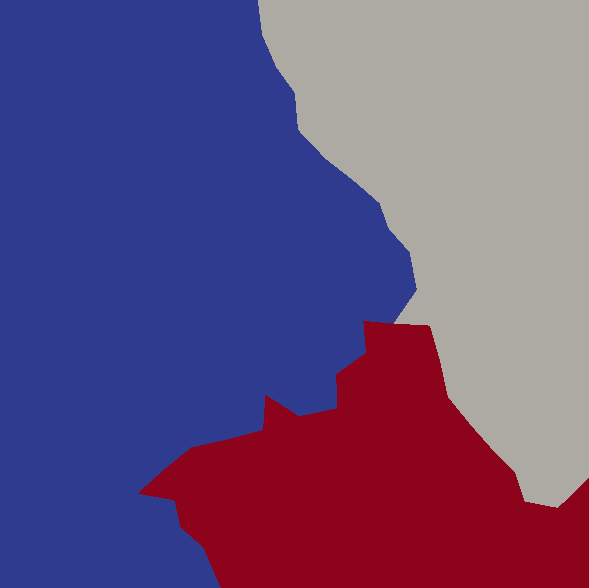}&\\
Kraus \cite{kraus_agglomeration-based_2004} & \includegraphics[height=2.5cm]{figures/simple2D_unstr_level_0} &\includegraphics[height=2.5cm]{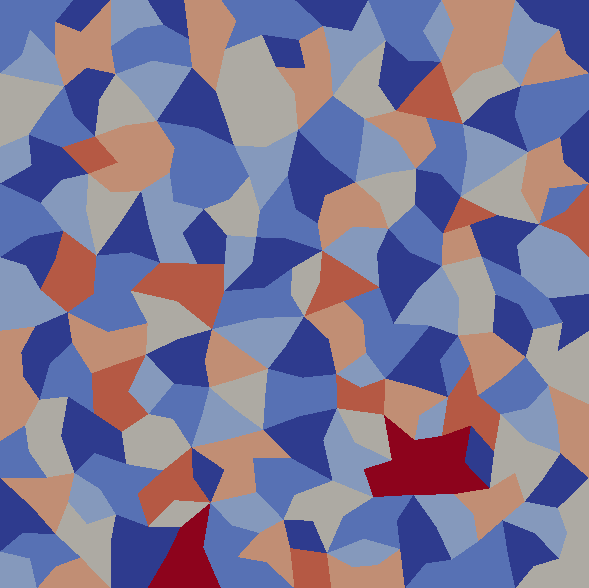}&\includegraphics[height=2.5cm]{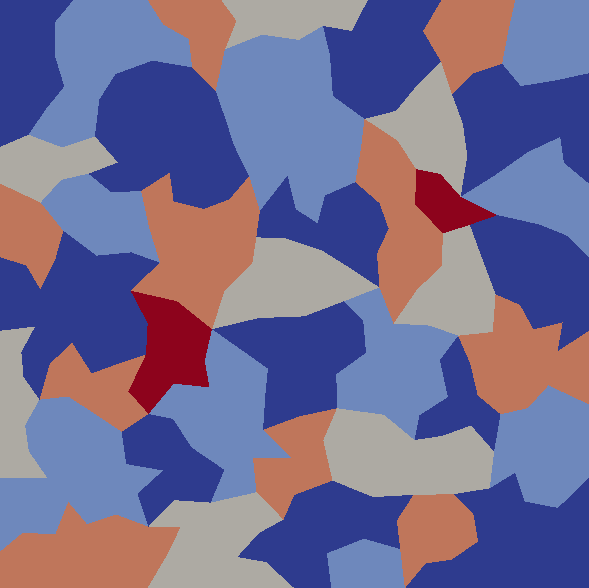}&\includegraphics[height=2.5cm]{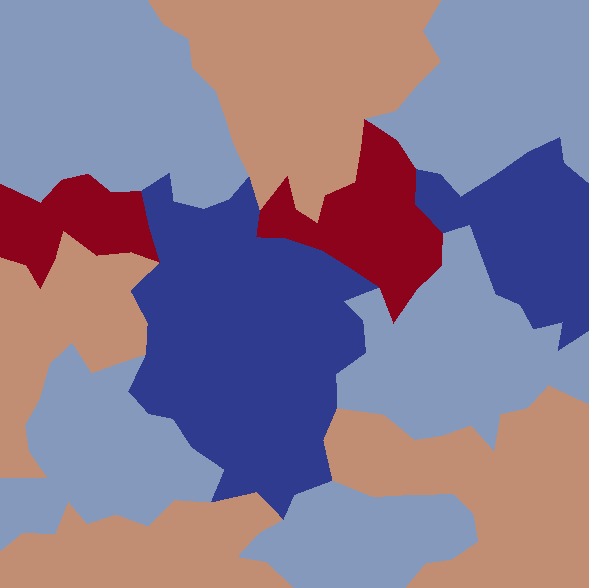}&\includegraphics[height=2.5cm]{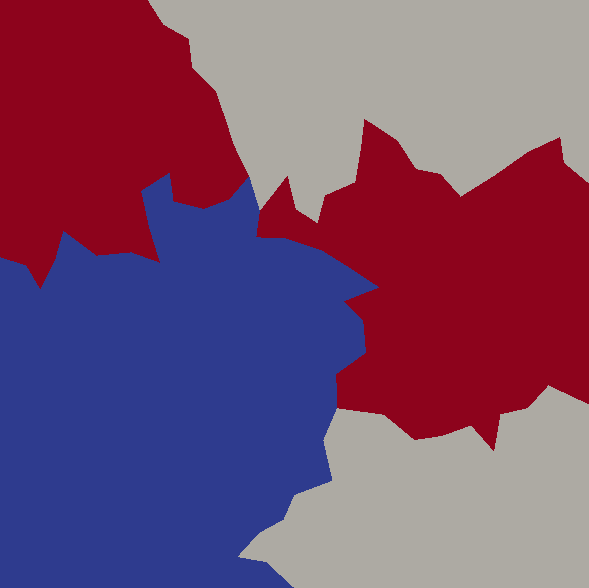}&\includegraphics[height=2.5cm]{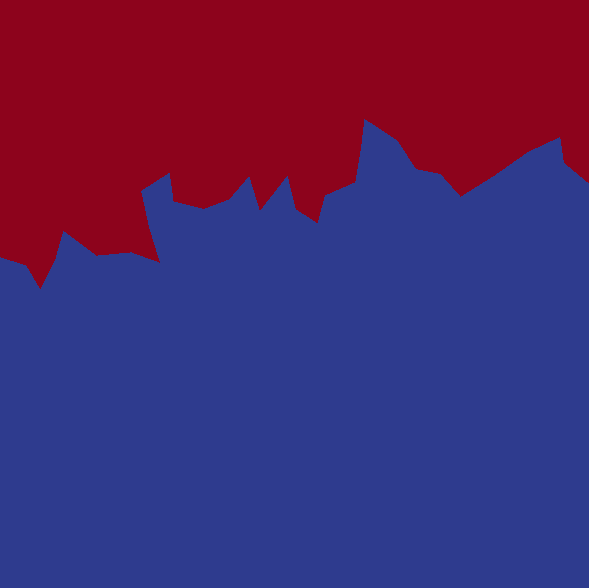}\\
RGB \cite{wabro_amgecoarsening_2006} & \includegraphics[height=2.5cm]{figures/simple2D_unstr_level_0} &\includegraphics[height=2.5cm]{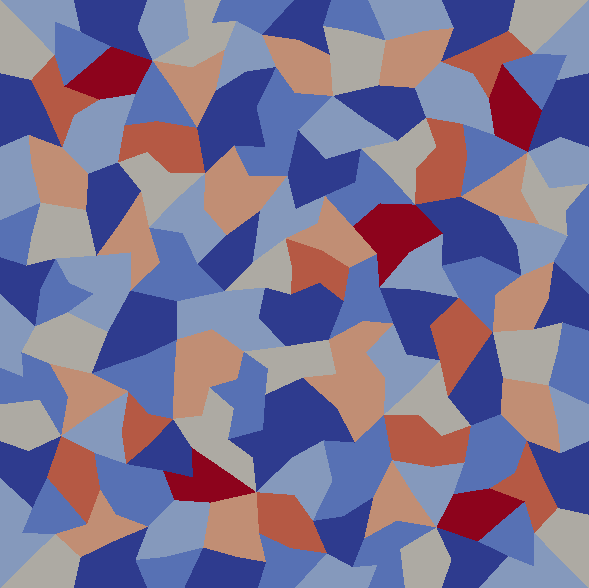}&\includegraphics[height=2.5cm]{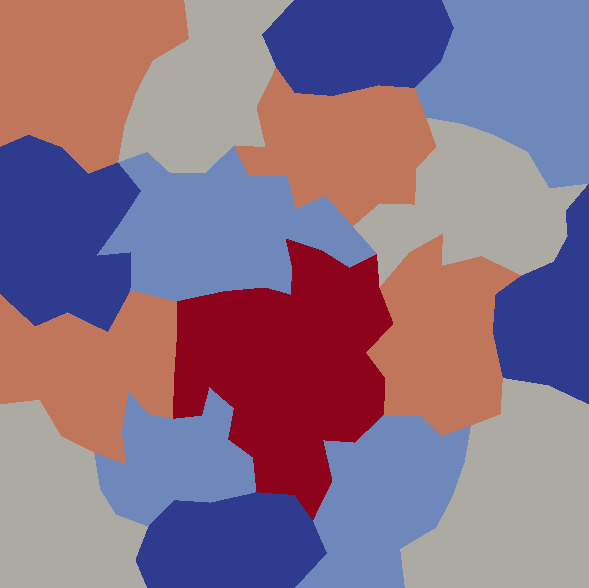}&\includegraphics[height=2.5cm]{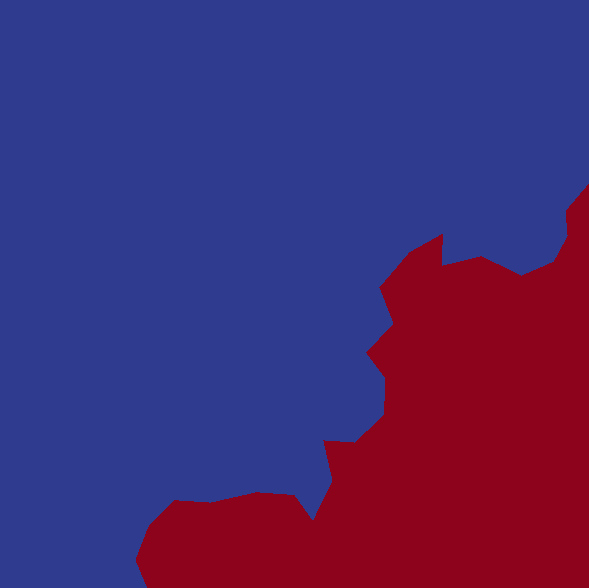}&&\\
Node & \includegraphics[height=2.5cm]{figures/simple2D_unstr_level_0} &\includegraphics[height=2.5cm]{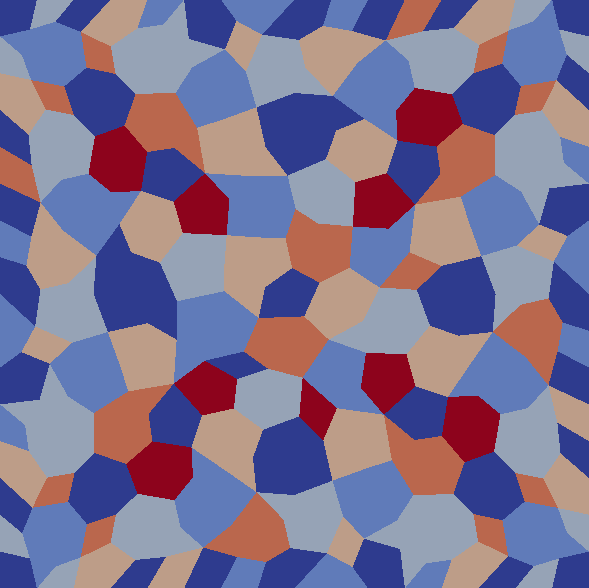}&\includegraphics[height=2.5cm]{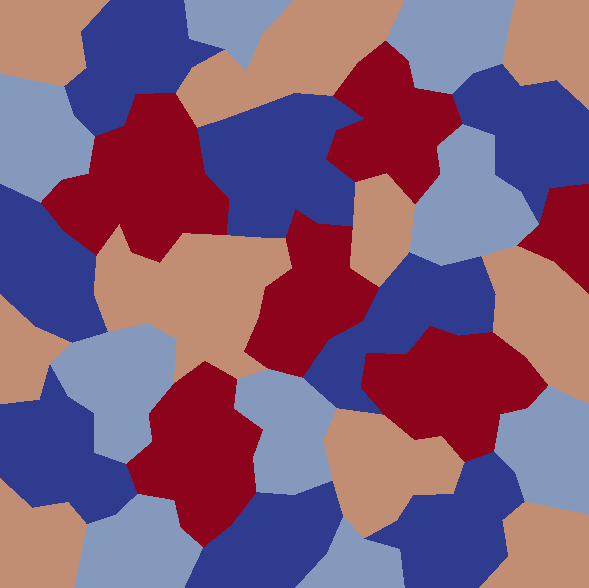}&\includegraphics[height=2.5cm]{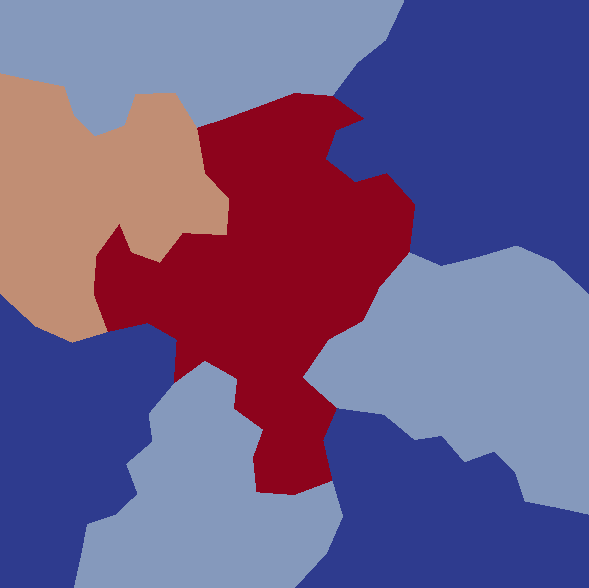}&&\\
Greedy & \includegraphics[height=2.5cm]{figures/simple2D_unstr_level_0} &\includegraphics[height=2.5cm]{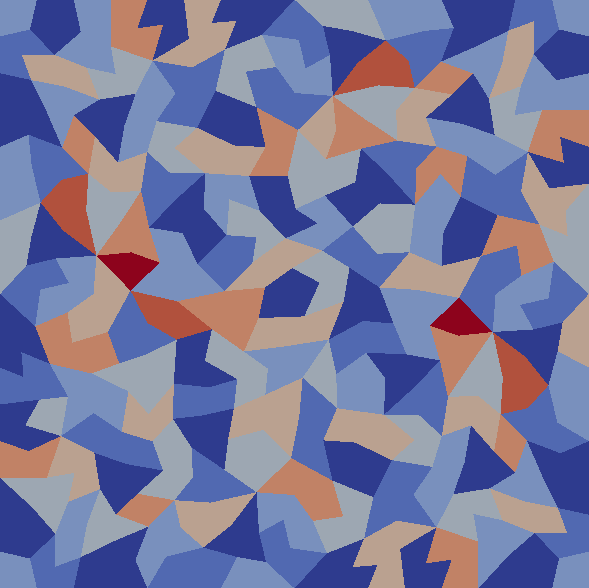}&\includegraphics[height=2.5cm]{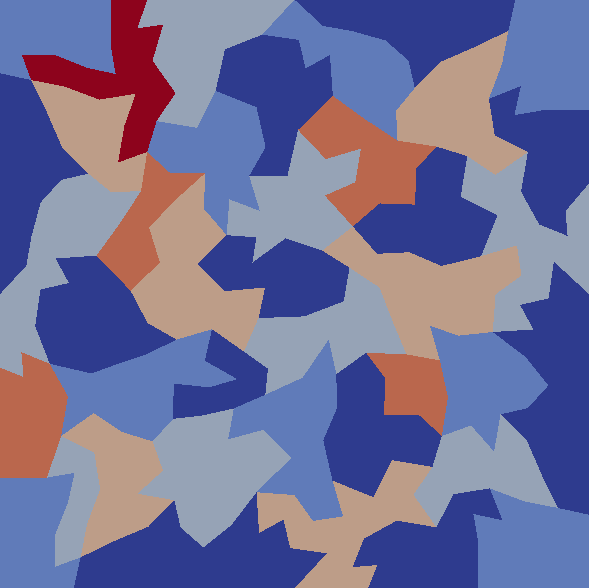}&\includegraphics[height=2.5cm]{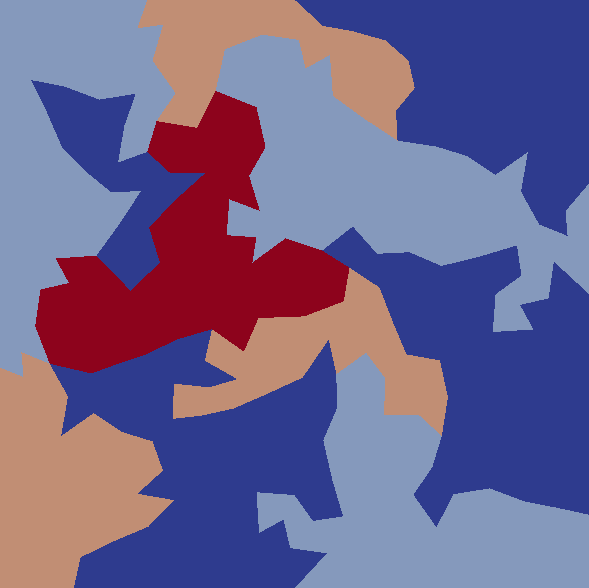}&\includegraphics[height=2.5cm]{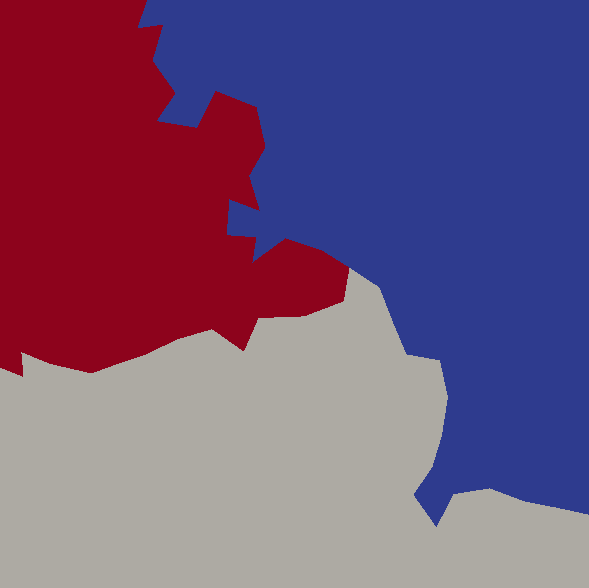}&\\
METIS \cite{karypis_metis-unstructured_1995} & \includegraphics[height=2.5cm]{figures/simple2D_unstr_level_0} &\includegraphics[height=2.5cm]{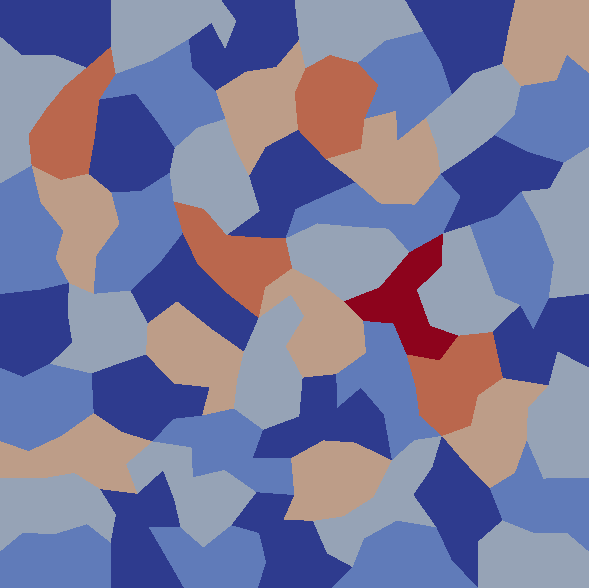}&\includegraphics[height=2.5cm]{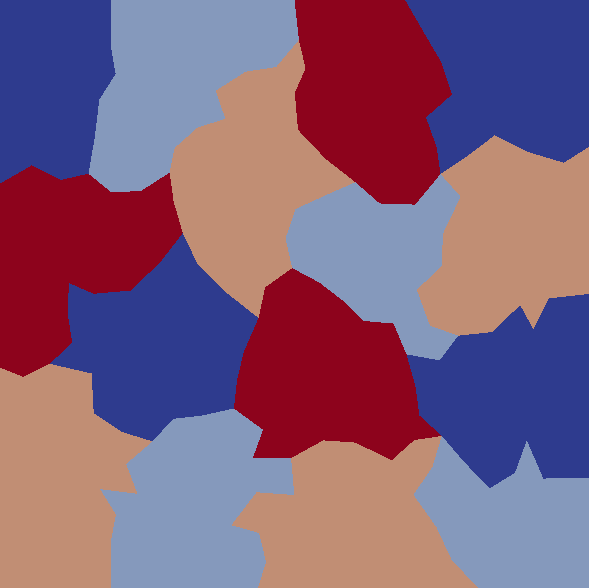}&\includegraphics[height=2.5cm]{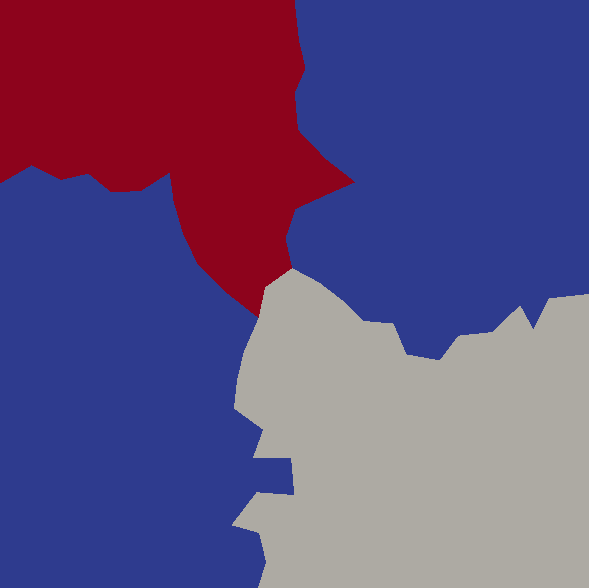}&&\\
MGridGen \cite{moulitsas_multilevel_2001} & \includegraphics[height=2.5cm]{figures/simple2D_unstr_level_0} &\includegraphics[height=2.5cm]{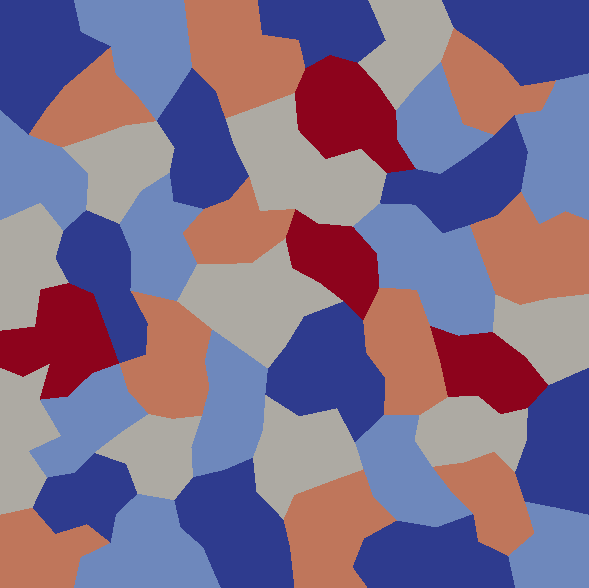}&\includegraphics[height=2.5cm]{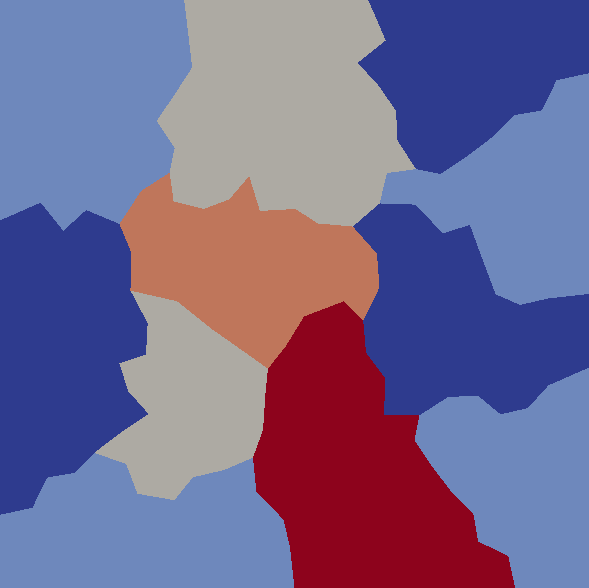}&\includegraphics[height=2.5cm]{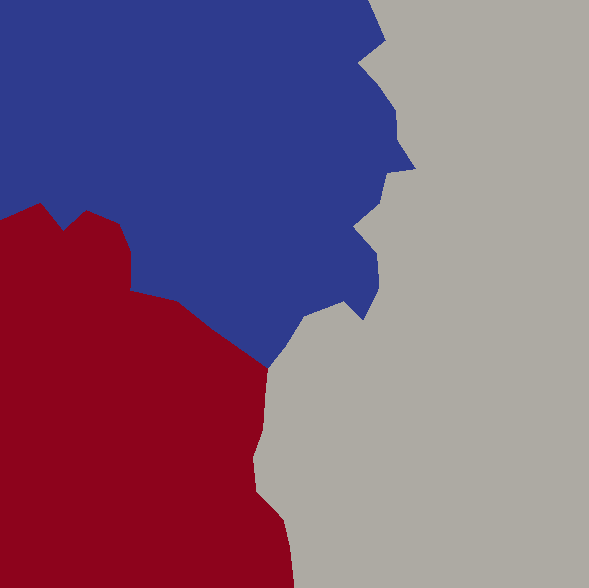}&&\\
\bottomrule
\end{tabular}
\end{center}
\caption{Different coarsening algorithms applied to a small (740 element) unstructured triangular mesh in 2D (n.b., this small mesh is not used in any of the example problems in this paper; it is merely for illustrative purposes)}
\label{tab:coarse_unstruct}
\end{table}

%\begin{table}[htp]
%\begin{center}
%\begin{tabular}{c c c c c}
%\toprule
%Alg. & Level 0 & 1 & 2 & 3\\
%\midrule
%1 & \includegraphics[height=2.5cm]{figures/simple3D_level_0} &\includegraphics[height=2.5cm]{figures/simple3D_m_1_level_1}&\includegraphics[height=2.5cm]{figures/simple3D_m_1_level_2}&\includegraphics[height=2.5cm]{figures/simple3D_m_1_level_3}\\
%2 & \includegraphics[height=2.5cm]{figures/simple3D_level_0} &\includegraphics[height=2.5cm]{figures/simple3D_m_2_level_1}&\includegraphics[height=2.5cm]{figures/simple3D_m_2_level_2}&\includegraphics[height=2.5cm]{figures/simple3D_m_2_level_3}\\
%3 & \includegraphics[height=2.5cm]{figures/simple3D_level_0} &\includegraphics[height=2.5cm]{figures/simple3D_m_3_level_1}&\includegraphics[height=2.5cm]{figures/simple3D_m_3_level_2}&\includegraphics[height=2.5cm]{figures/simple3D_m_3_level_3}\\
%4 & \includegraphics[height=2.5cm]{figures/simple3D_level_0} &\includegraphics[height=2.5cm]{figures/simple3D_m_4_level_1}&\includegraphics[height=2.5cm]{figures/simple3D_m_4_level_2}&\\
%5 & \includegraphics[height=2.5cm]{figures/simple3D_level_0} &\includegraphics[height=2.5cm]{figures/simple3D_m_5_level_1}&\includegraphics[height=2.5cm]{figures/simple3D_m_5_level_2}&\\
%6 & \includegraphics[height=2.5cm]{figures/simple3D_level_0} &\includegraphics[height=2.5cm]{figures/simple3D_m_6_level_1}&\includegraphics[height=2.5cm]{figures/simple3D_m_6_level_2}&\\
%\bottomrule
%\end{tabular}
%\end{center}
%\label{fig:coarse_struct3D}
%\caption{Coarsening algorithms applied to a small (512 element) structured grid in 3D}
%\end{table}

\begin{table}[p]
\begin{center}
\begin{tabular}{>{\centering\arraybackslash} m{1.75cm} >{\centering\arraybackslash} m{2.5cm} >{\centering\arraybackslash} m{2.5cm} >{\centering\arraybackslash} m{2.5cm} >{\centering\arraybackslash} m{2.5cm} >{\centering\arraybackslash} m{2.5cm} >{\centering\arraybackslash} m{2.5cm}}
\toprule
Alg. & Orig. mesh & Level 1 & 2 & 3 & 4 & 5 \\
\midrule
Jones \cite{jones_amge_2001} & \includegraphics[height=2.5cm]{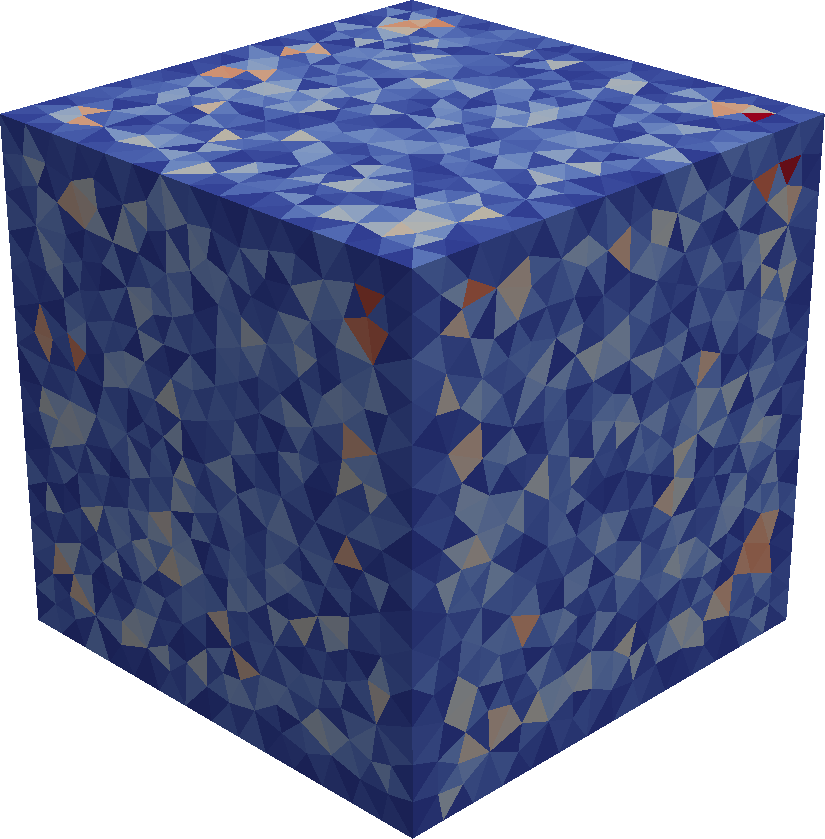} &\includegraphics[height=2.5cm]{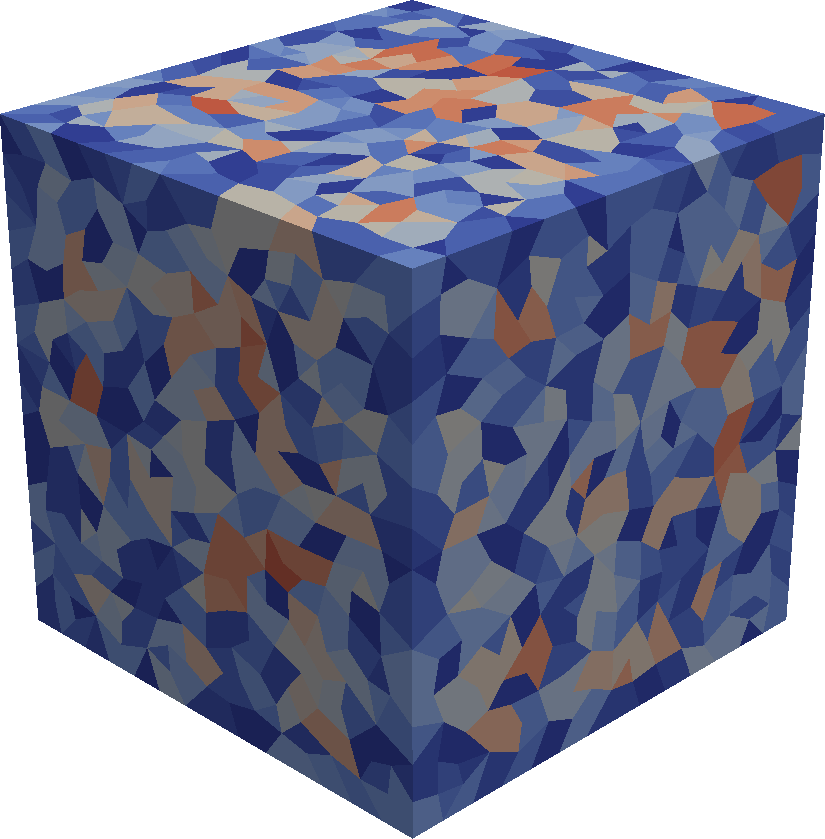}&\includegraphics[height=2.5cm]{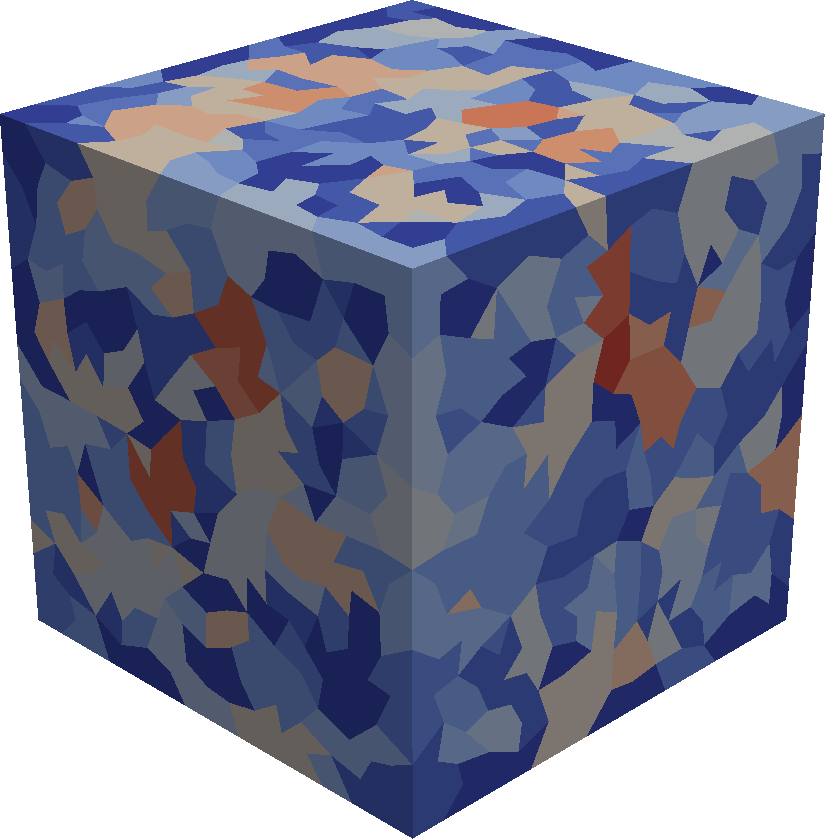}&\includegraphics[height=2.5cm]{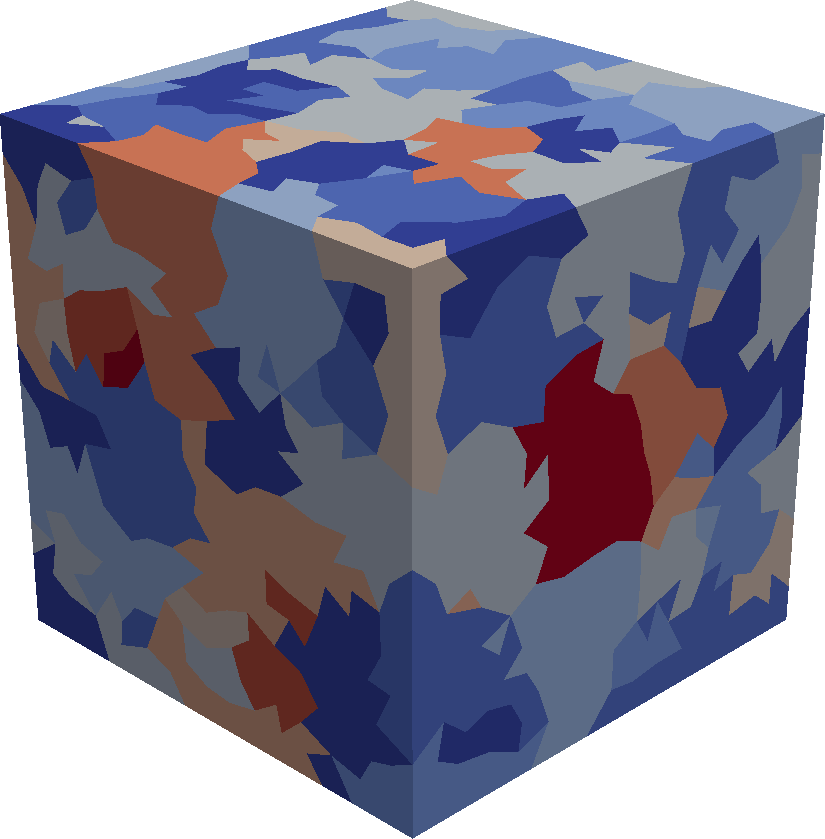}&\includegraphics[height=2.5cm]{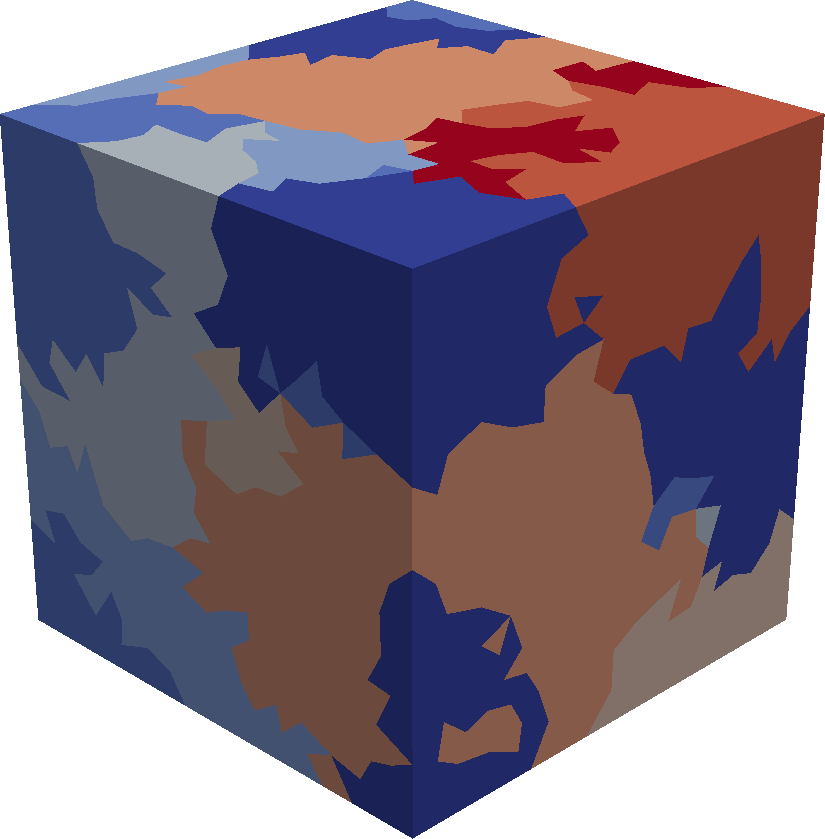}&\includegraphics[height=2.5cm]{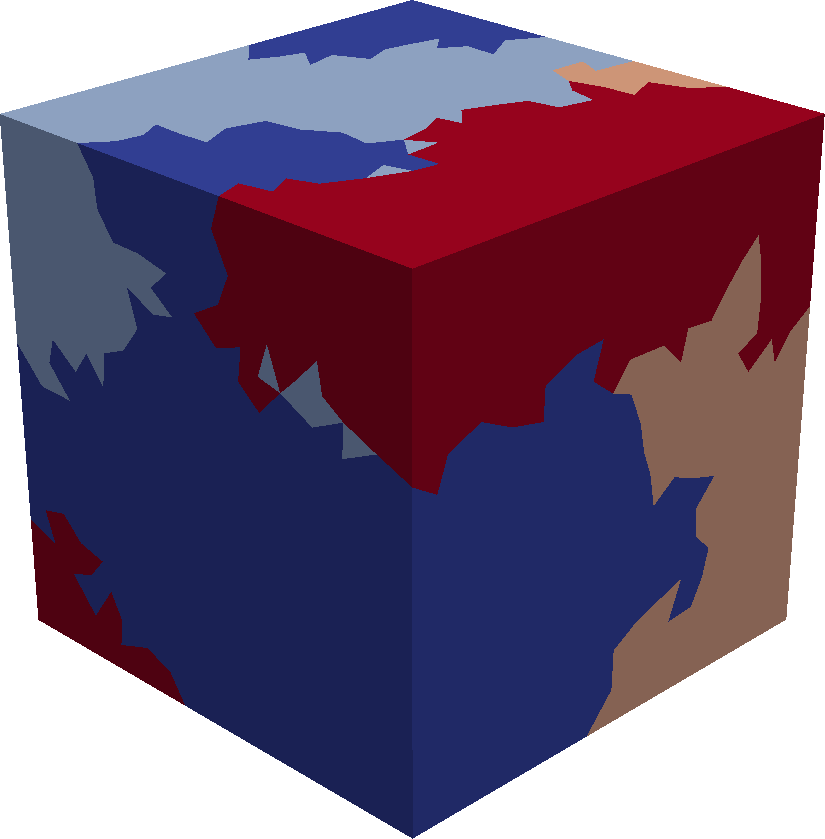}\\
Kraus \cite{kraus_agglomeration-based_2004} & \includegraphics[height=2.5cm]{figures/c_0} &\includegraphics[height=2.5cm]{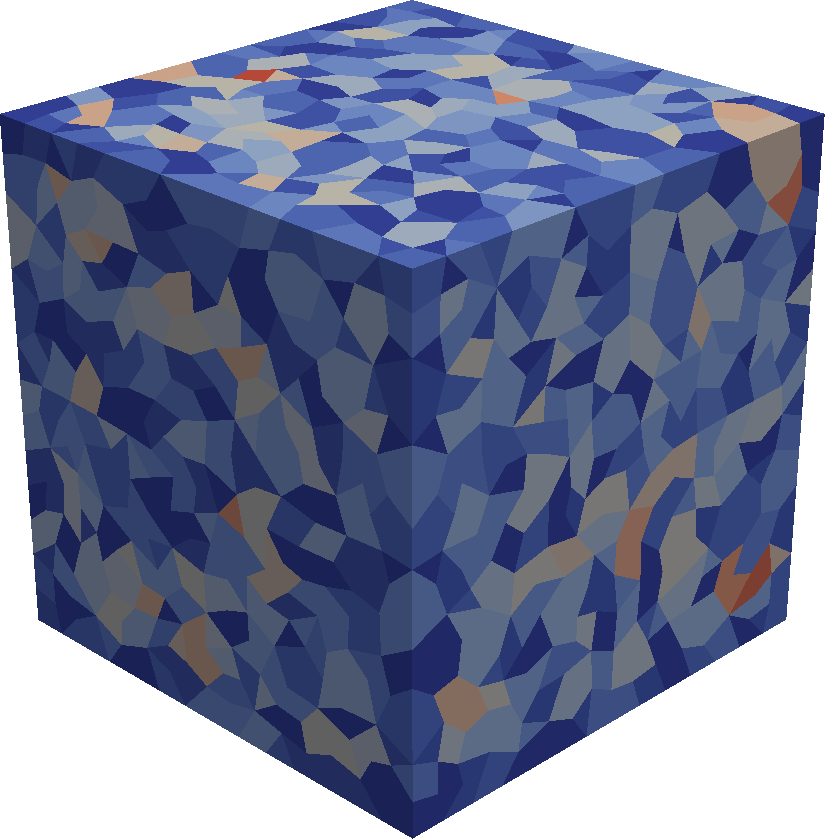}&\includegraphics[height=2.5cm]{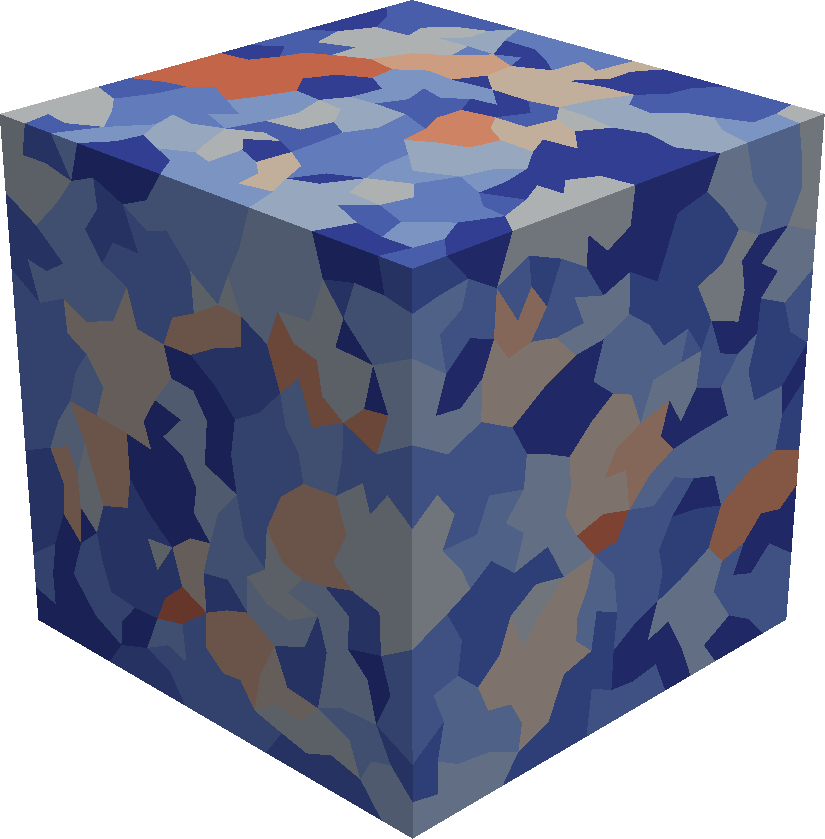}&\includegraphics[height=2.5cm]{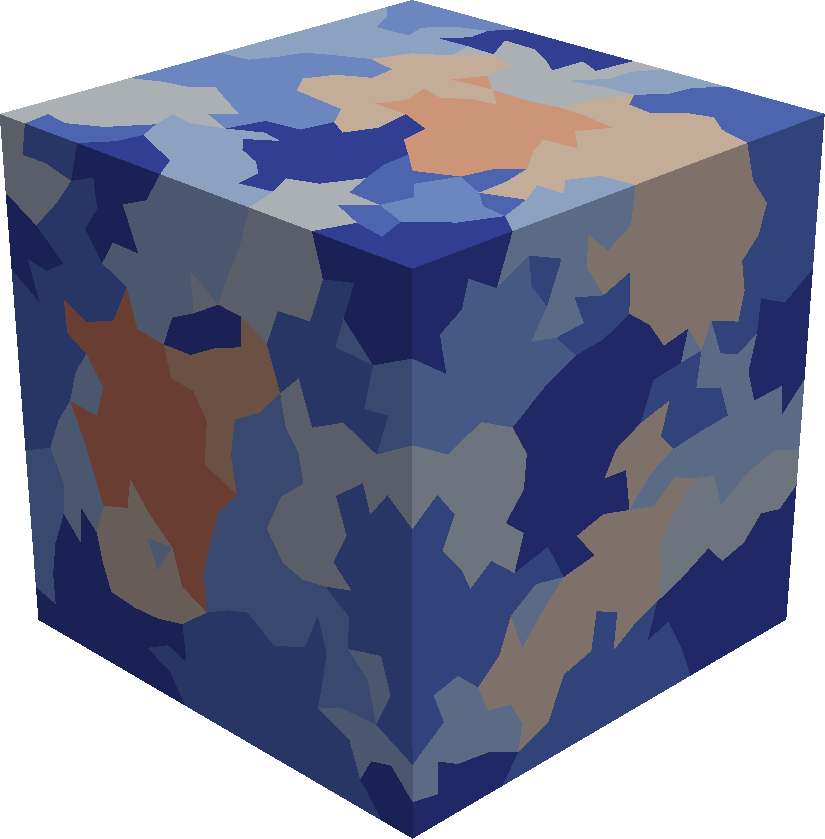}&\includegraphics[height=2.5cm]{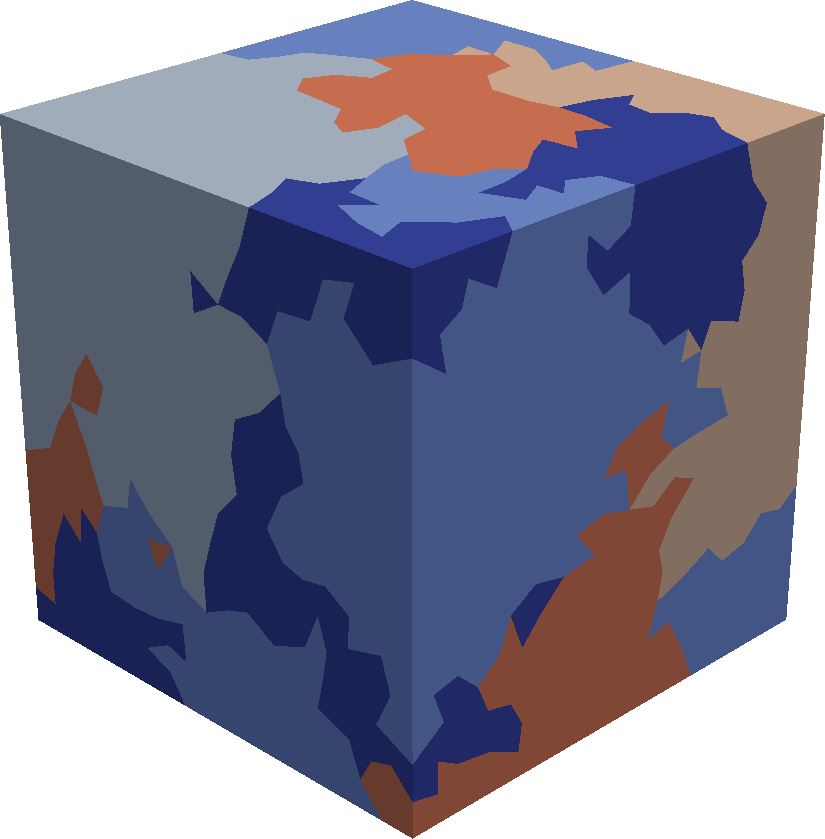}&\includegraphics[height=2.5cm]{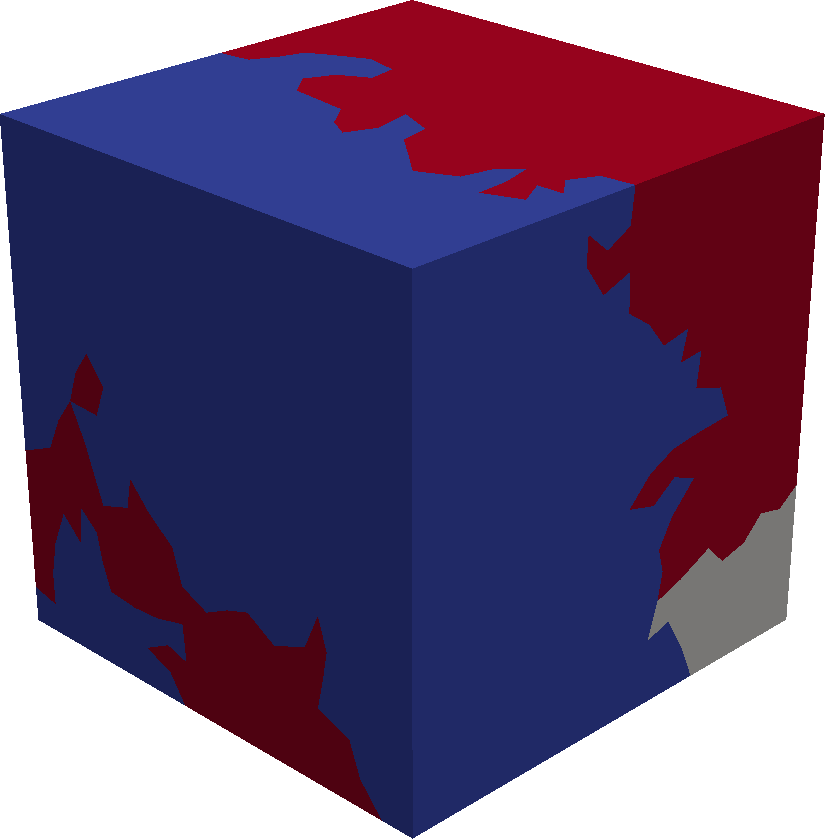}\\
RGB \cite{wabro_amgecoarsening_2006} & \includegraphics[height=2.5cm]{figures/c_0} &\includegraphics[height=2.5cm]{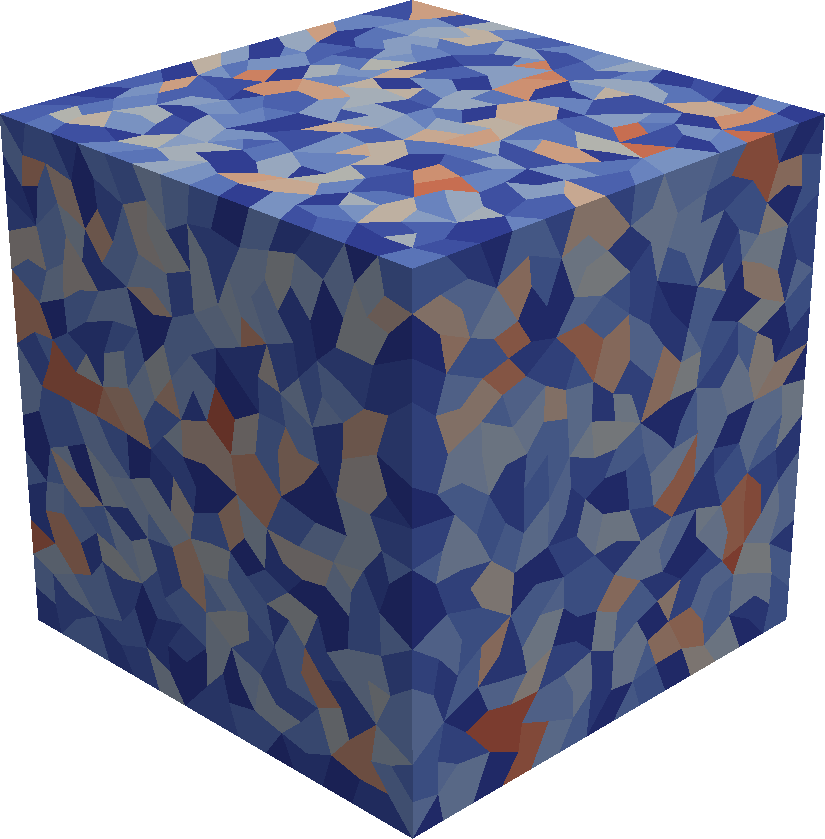}&\includegraphics[height=2.5cm]{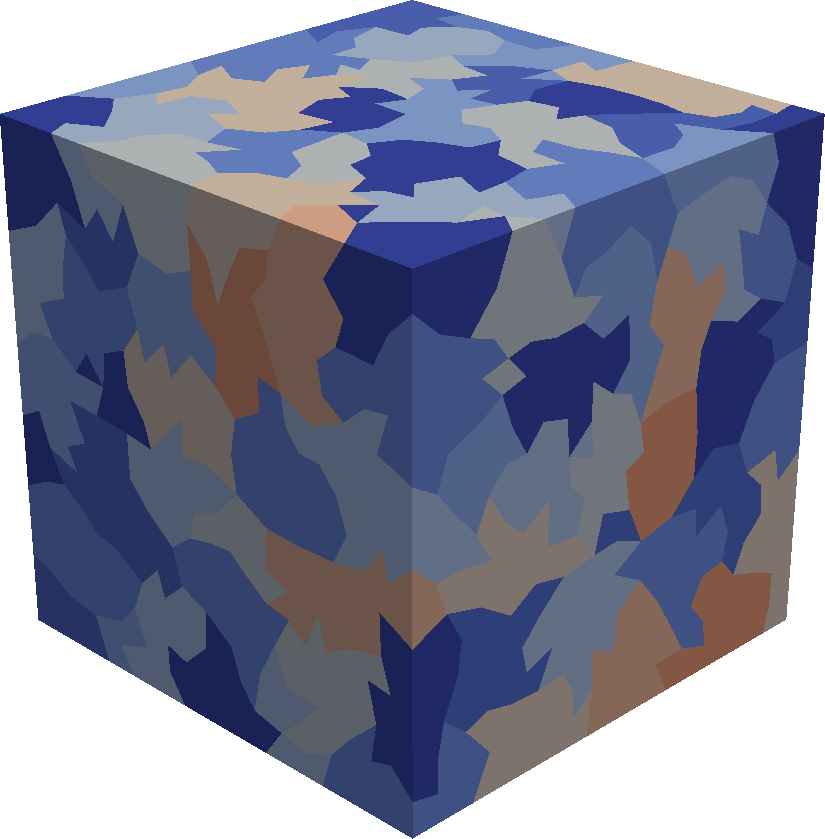}&\includegraphics[height=2.5cm]{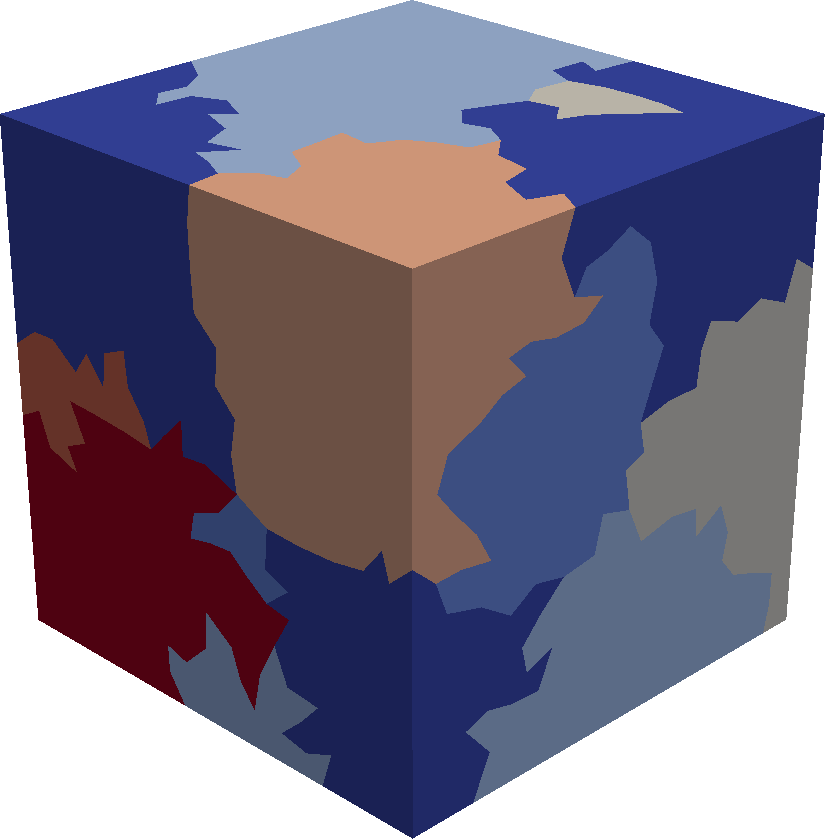}&&\\
Node & \includegraphics[height=2.5cm]{figures/c_0} & \includegraphics[height=2.5cm]{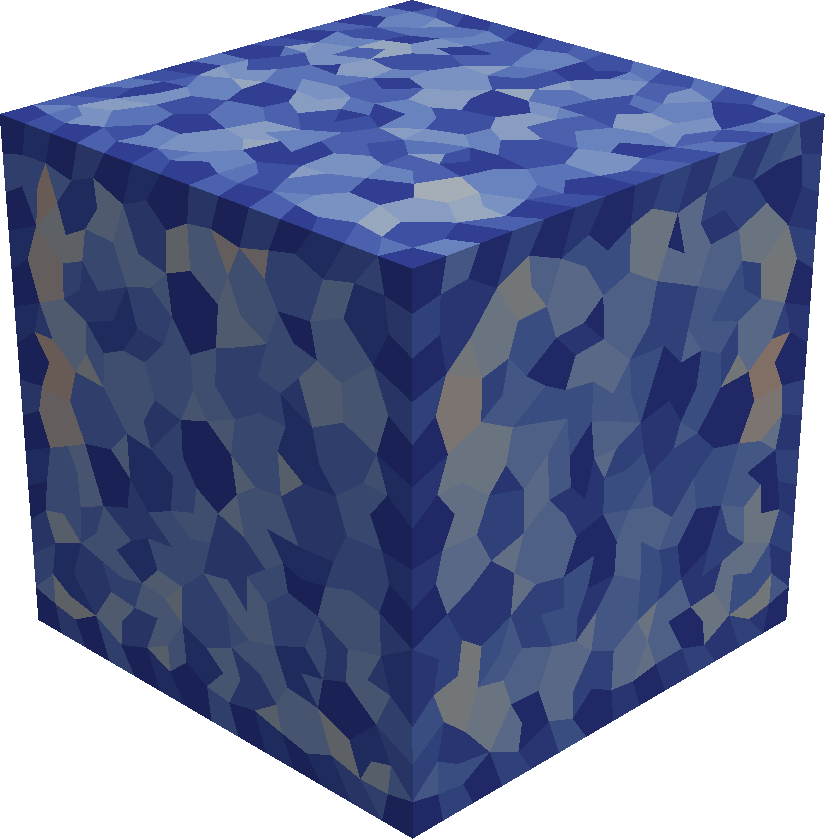}&\includegraphics[height=2.5cm]{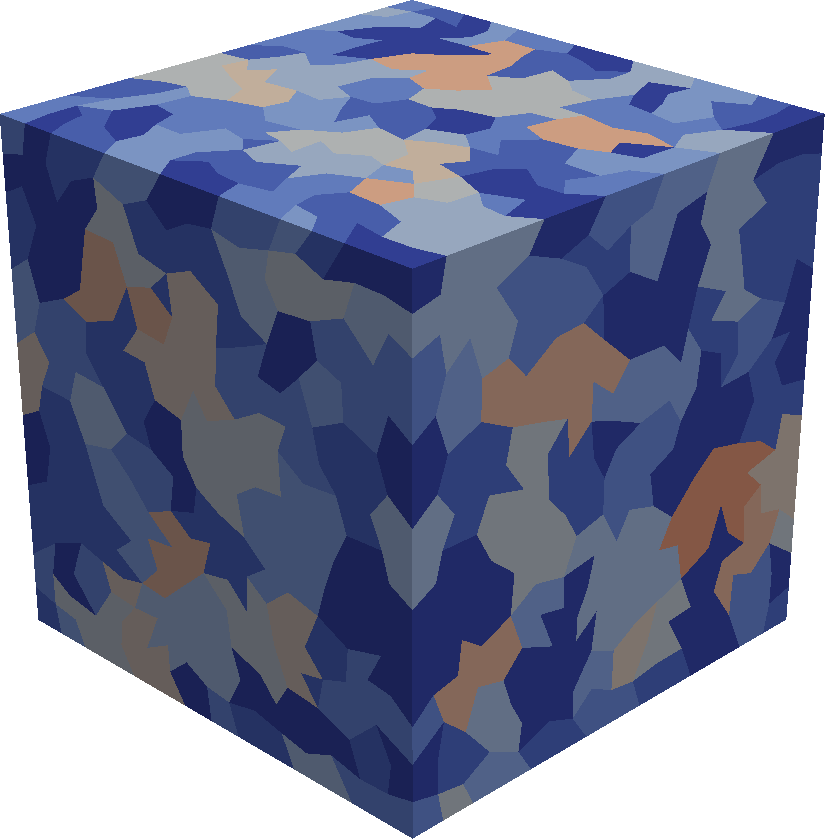}&\includegraphics[height=2.5cm]{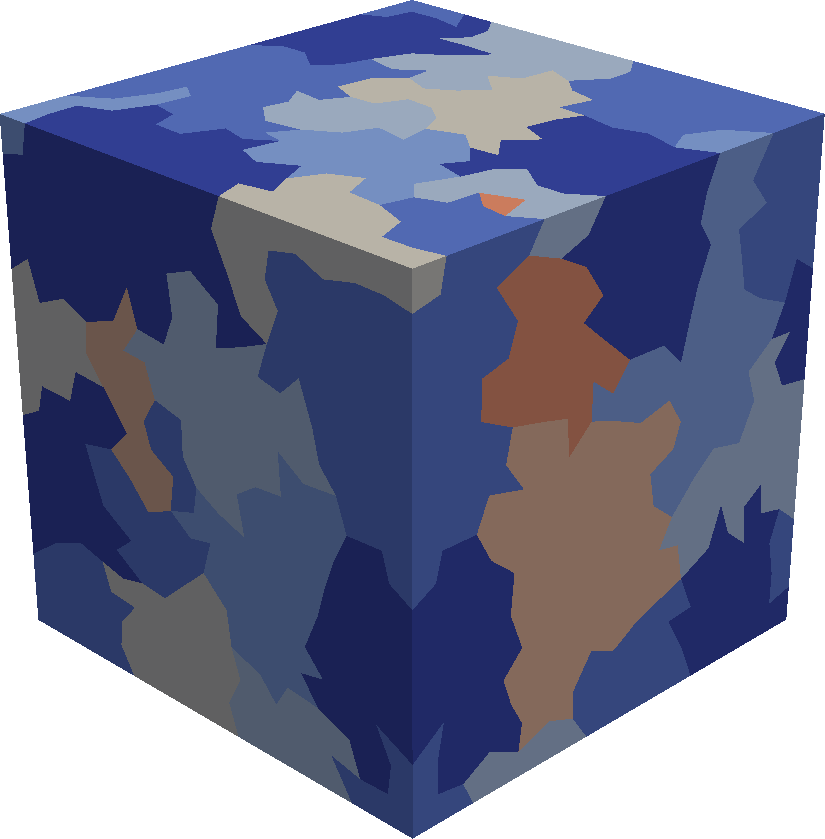}&\includegraphics[height=2.5cm]{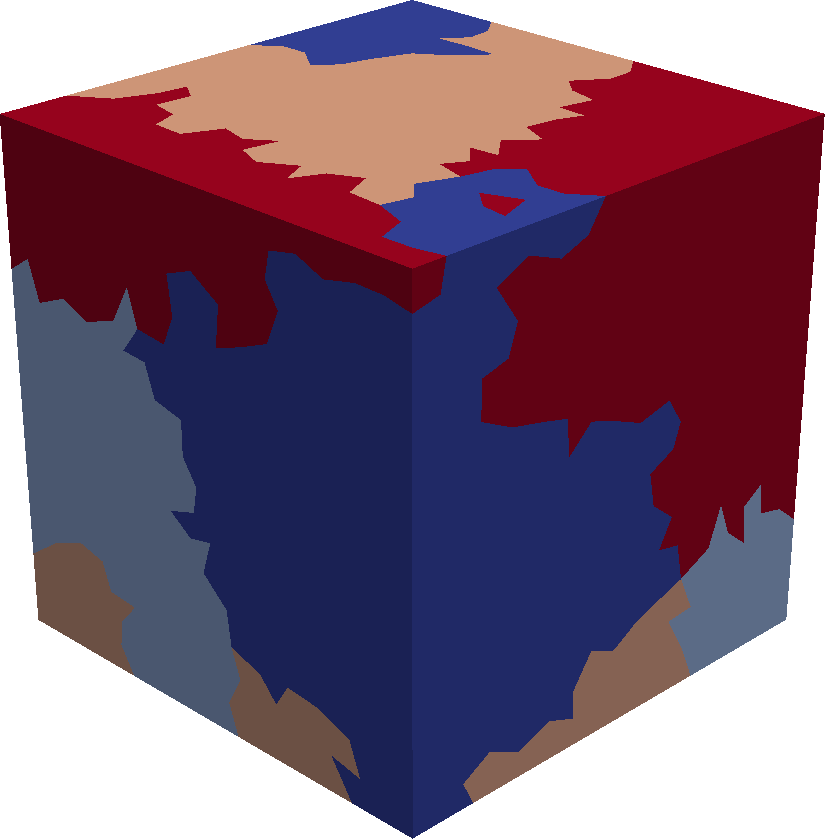}&\\
Greedy & \includegraphics[height=2.5cm]{figures/c_0} &\includegraphics[height=2.5cm]{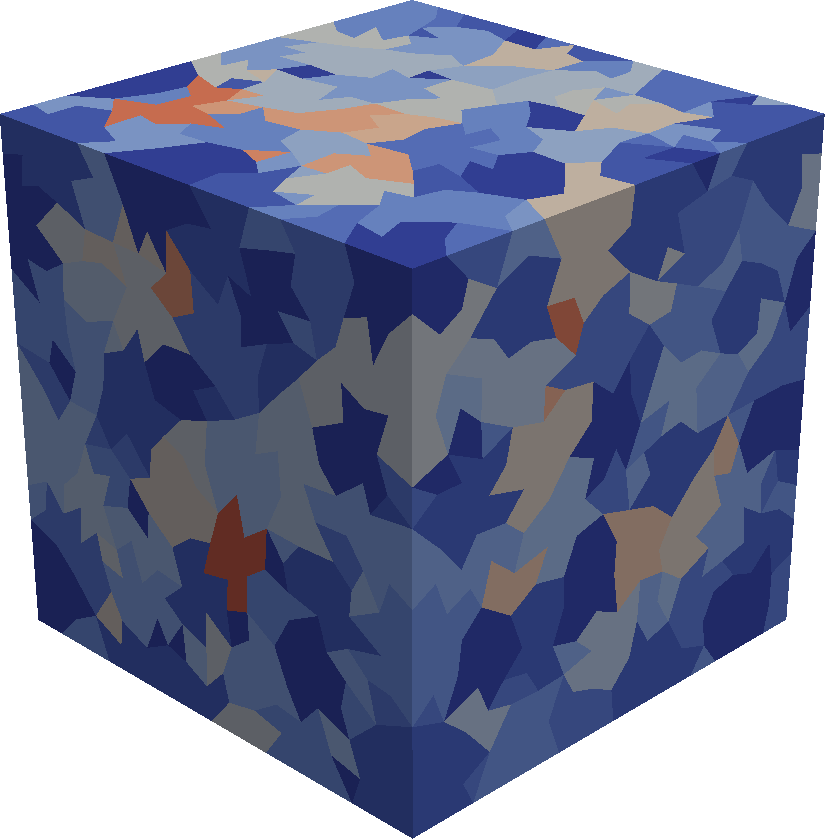}&\includegraphics[height=2.5cm]{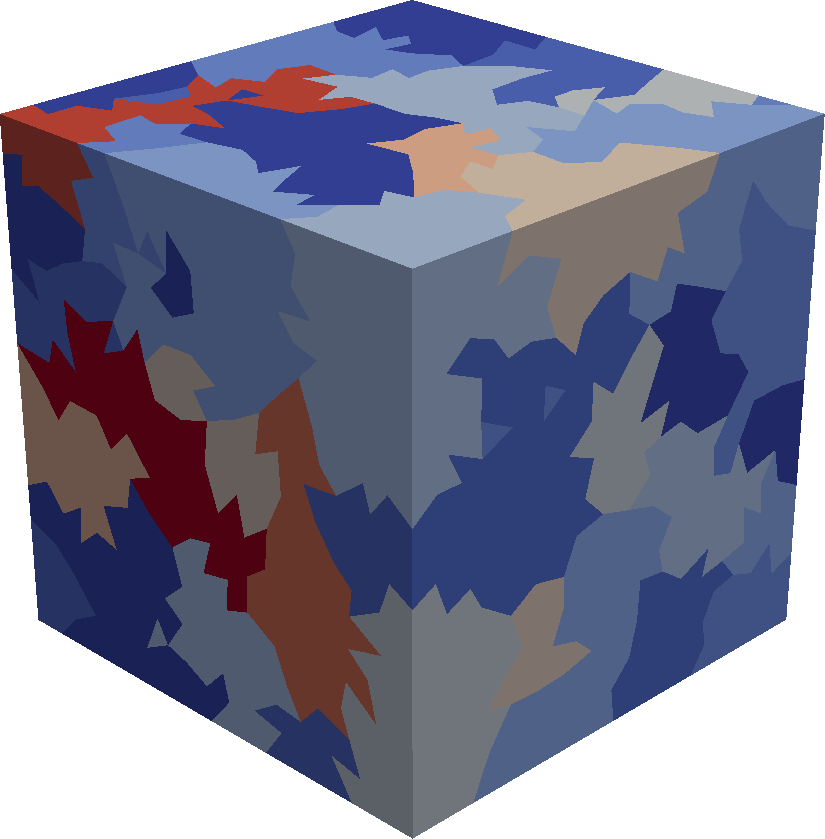}&\includegraphics[height=2.5cm]{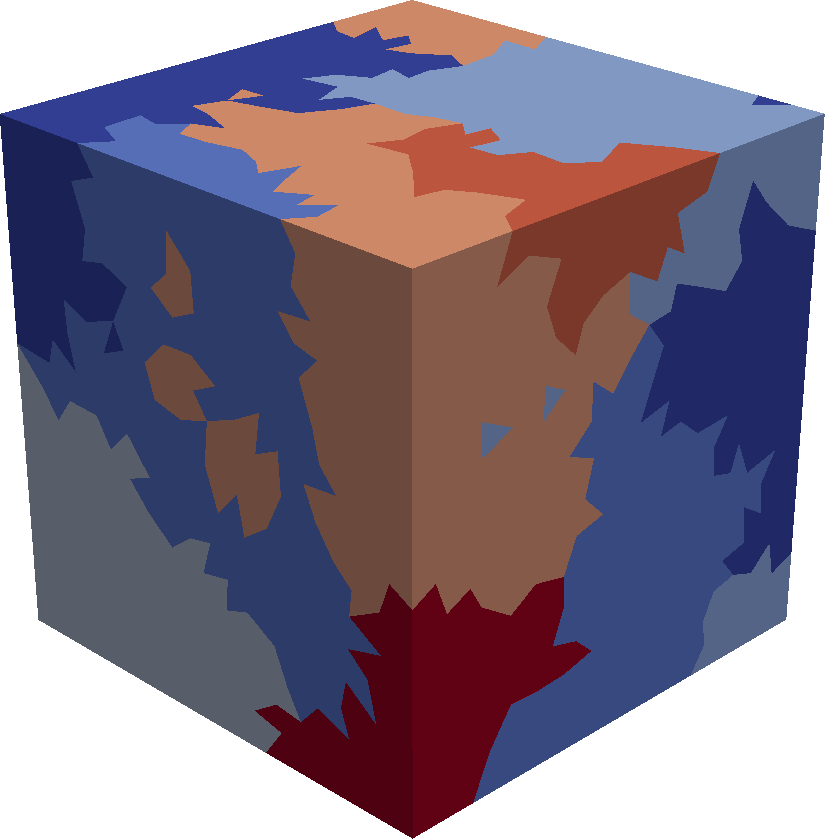}&\includegraphics[height=2.5cm]{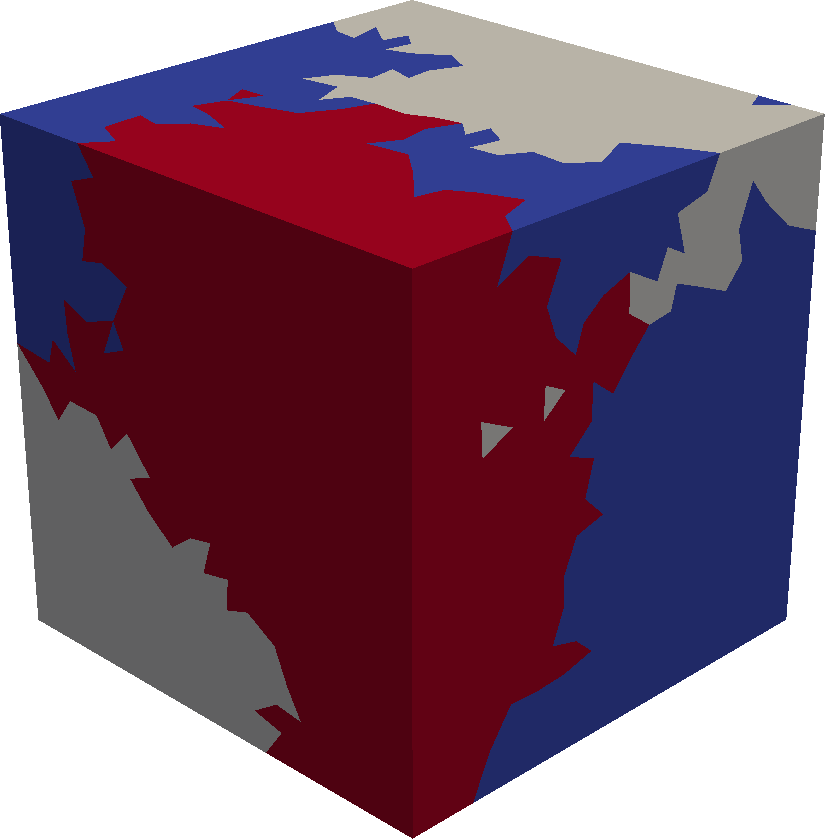}&\\
METIS \cite{karypis_metis-unstructured_1995} & \includegraphics[height=2.5cm]{figures/c_0} &\includegraphics[height=2.5cm]{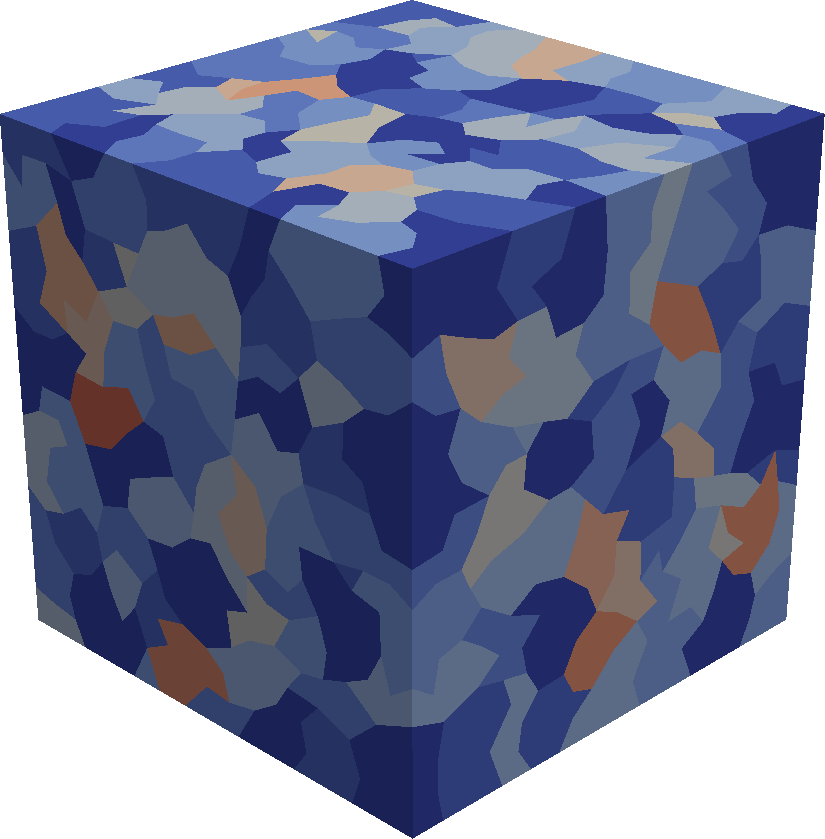}&\includegraphics[height=2.5cm]{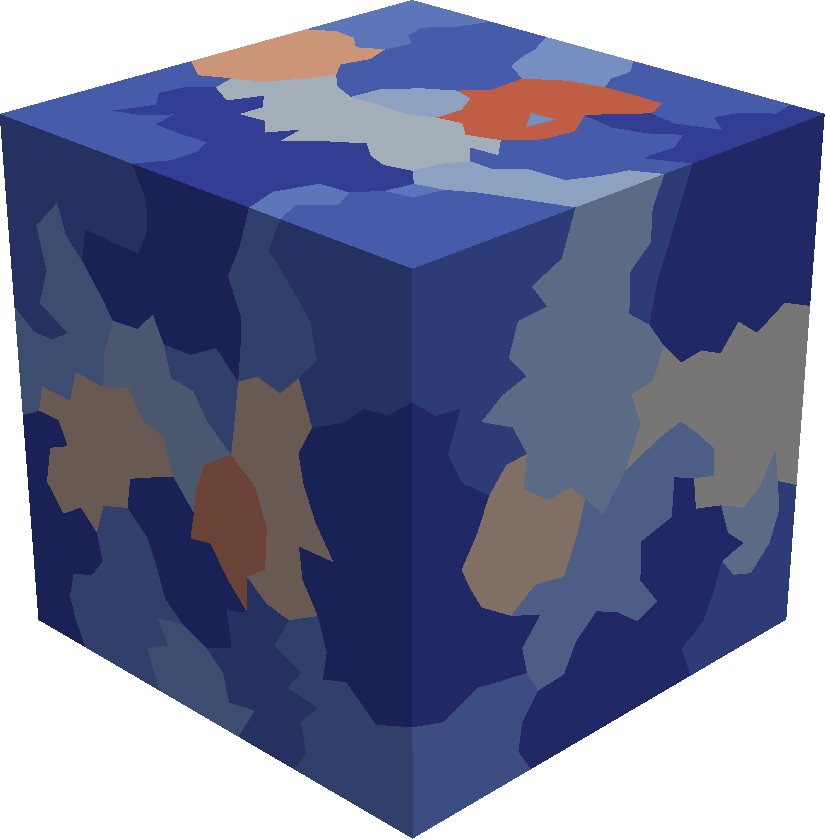}&\includegraphics[height=2.5cm]{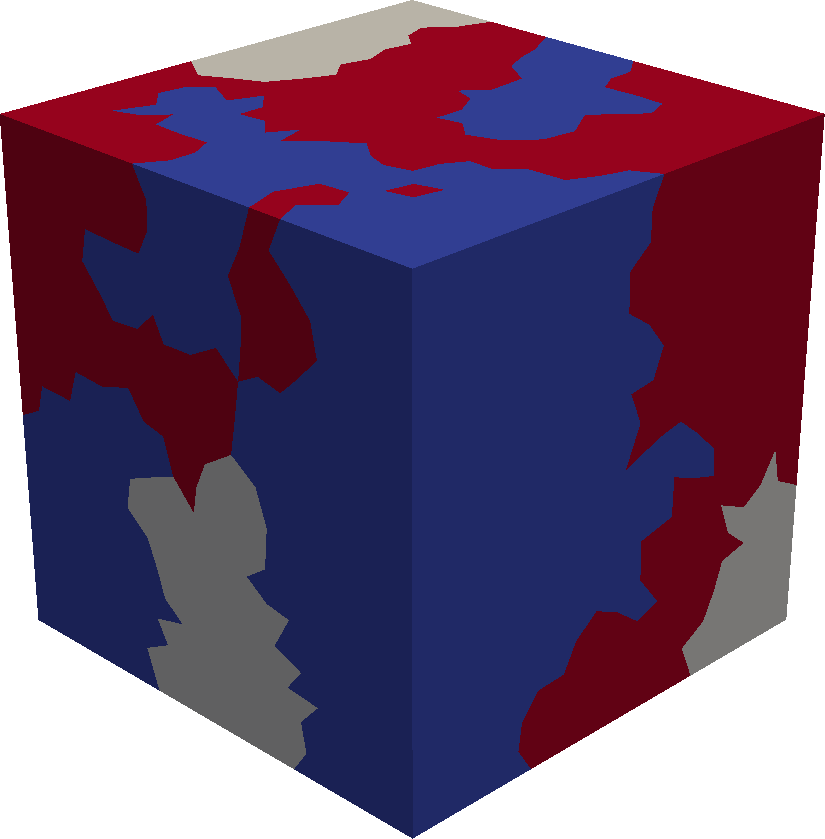}&&\\
MGridGen \cite{moulitsas_multilevel_2001} & \includegraphics[height=2.5cm]{figures/c_0} &\includegraphics[height=2.5cm]{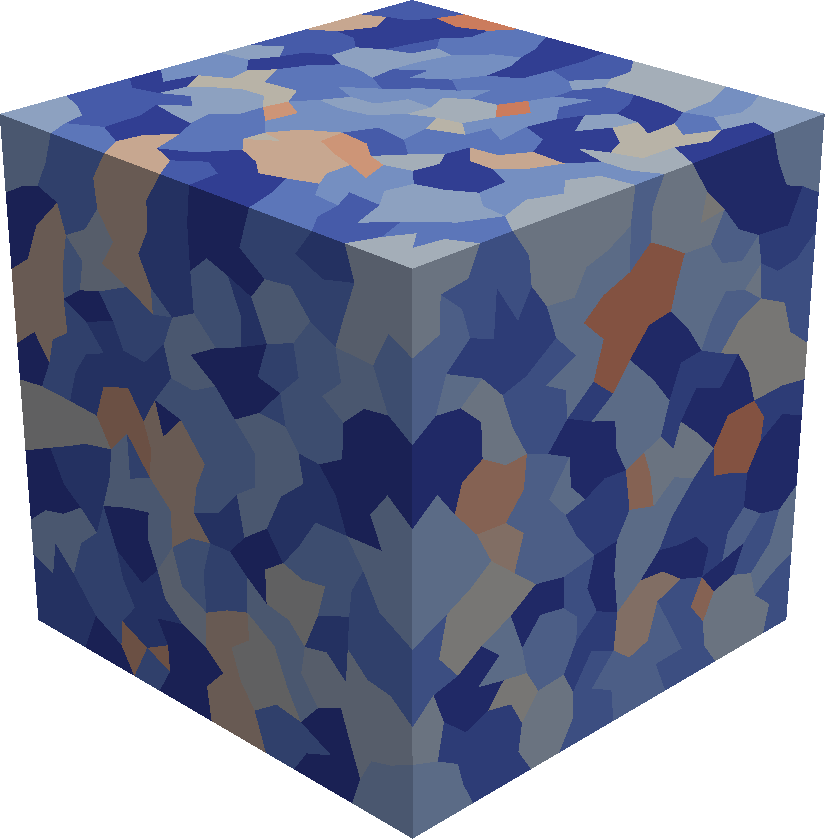}&\includegraphics[height=2.5cm]{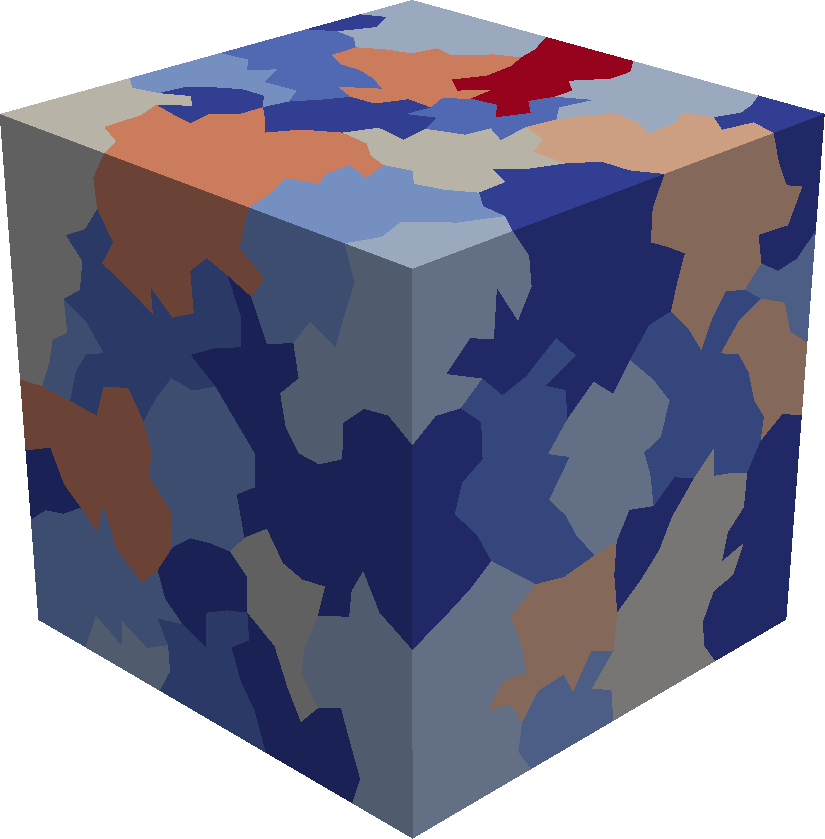}&\includegraphics[height=2.5cm]{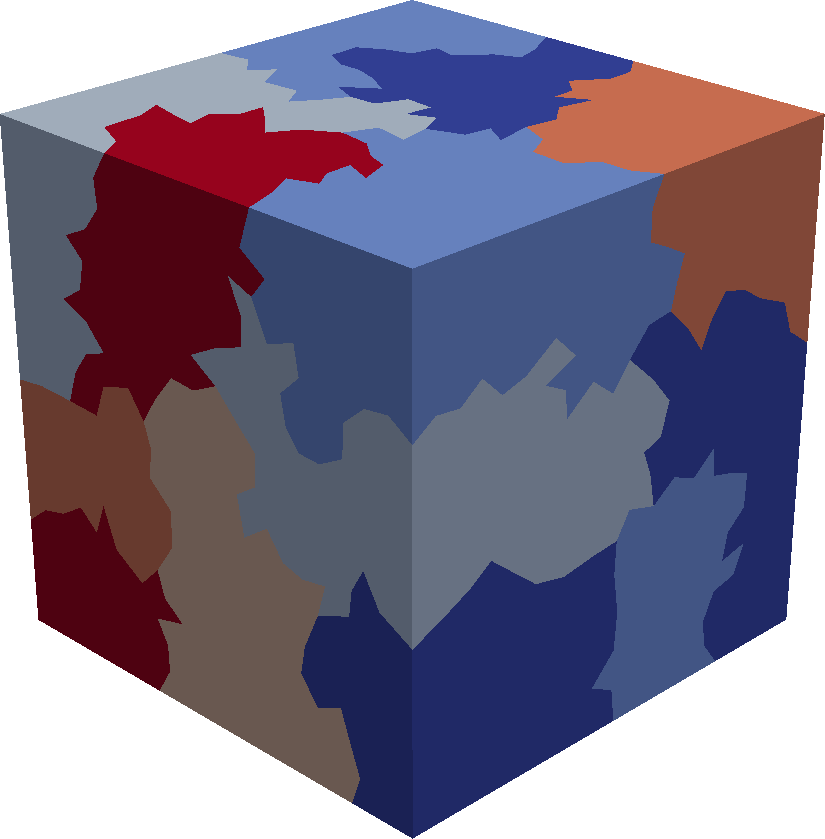}&\includegraphics[height=2.5cm]{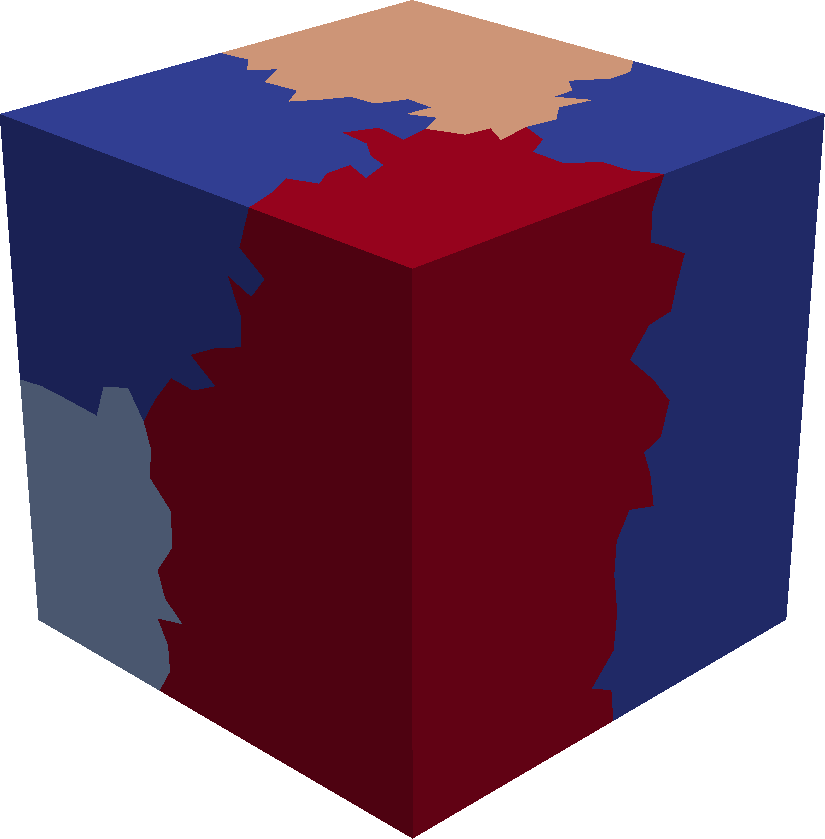}&\\
\bottomrule
\end{tabular}
\end{center}
\caption{Different coarsening algorithms applied to a small (12876 element) unstructured tetrahedron-based grid in 3D (n.b., this small mesh is not used in any of the example problems in this paper; it is merely for illustrative purposes)}
\label{tab:coarse_unstruct3D}
\end{table}
%-----------------------------
\section{Results}
\label{sec:Results}
The results shown in this section focus entirely on using agglomeration algorithms on unstructured meshes. We use gmsh \cite{geuzaine_gmsh_2009} to mesh all of our problem domains, using triangles in 2D and tetrahedrons in 3D. We have not investigated the use of different mesh generators in this work, though in a general sense, the results shown below tend to rely on simple trends that hold on all unstructured triangular grids (like increased connectivity). 

In all the results shown below, the solve time (s) refers only to the time spent in by our iterative method to solve to a relative or absolute tolerance of 10$^{-10}$. We use an initial guess of zero in our iterative method. The multigrid setup time (s) refers to the time spent performing the element coarsening algorithms, any cleanup, assembly of coarse topologies (e.g., coarse faces or coarse edges) and the projection of spatial tables onto the lower grids. This time does not include the projection of the discretised form of \eref{eq:bte} onto the lower grids given that is not needed with our matrix-free approach, and we should note it also does not include the projection of our stabilisation matrix $\mat{D}^{-1}$ onto the lower grids (see \cite{buchan_inner-element_2010, Dargaville2019}; the cost of this is proportional to the cost of projecting our spatial tables). Finally the runtime (s) refers to the total runtime of our code, including all the steps mentioned above and any other miscellaneous time. 

In this section, we use two related metrics to quantify the properties of our unstructured meshes. It is useful to consider these in terms of an optimal coarsening on a structured quadrilateral mesh (ignoring boundary effects). These are:
\begin{enumerate}
\item The average connectivity on a mesh, which is the average number of elements each node on the mesh is connected to. For a structured quadrilateral mesh in 2D, the average connectivity is 4 (8 in 3D), and this remains constant on lower grids with an optimal coarsening.
\item The node/element ratio on a set level. For a structured quadrilateral mesh in 2D, this remains constant at 1 with an optimal coarsening. 
\end{enumerate}  
%--------------
\subsection{Choosing ideal agglomerate size}
\label{sec:Choosing}
Three of the algorithms considered in \secref{sec:Coarsening algorithms} (Greedy, METIS and MGridGen) require the input of a desired agglomerate size (we denote these algorithms as ``size-based''). Before comparing the results between all 7 coarsening algorithms, we must first attempt to determine an ideal agglomerate size for these size-based algorithms when used in a geometric context, like our work.

On a structured quadrilateral mesh, the ideal agglomerate size is easy to determine and constant. On unstructured grids, the increased connectivity makes the choice of ideal agglomerate size difficult and there is no reason that this ideal size should remain the same on different levels. In order to try and quantify this size, we use the size-based coarsening algorithms on two simple test problems in both 2 and 3 dimensions, while varying the desired agglomerate size.

\begin{figure}[th]
\centering
\subfloat[][2D - red region is source]{\label{fig:2D_eg}\includegraphics[width =0.425\textwidth]{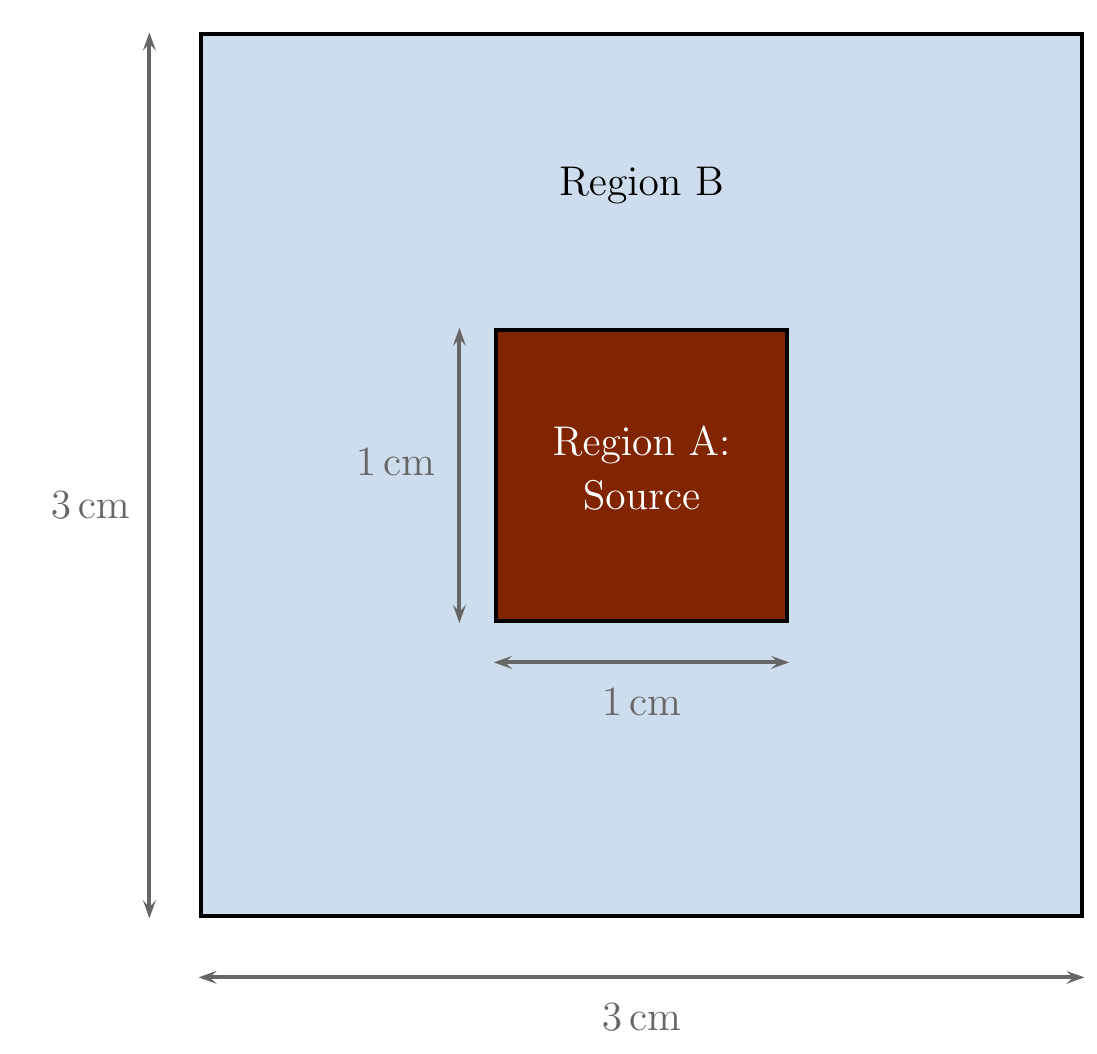}} \hspace{0.5cm}
\subfloat[][3D - red region is source]{\label{fig:3D_eg}\includegraphics[width =0.425\textwidth]{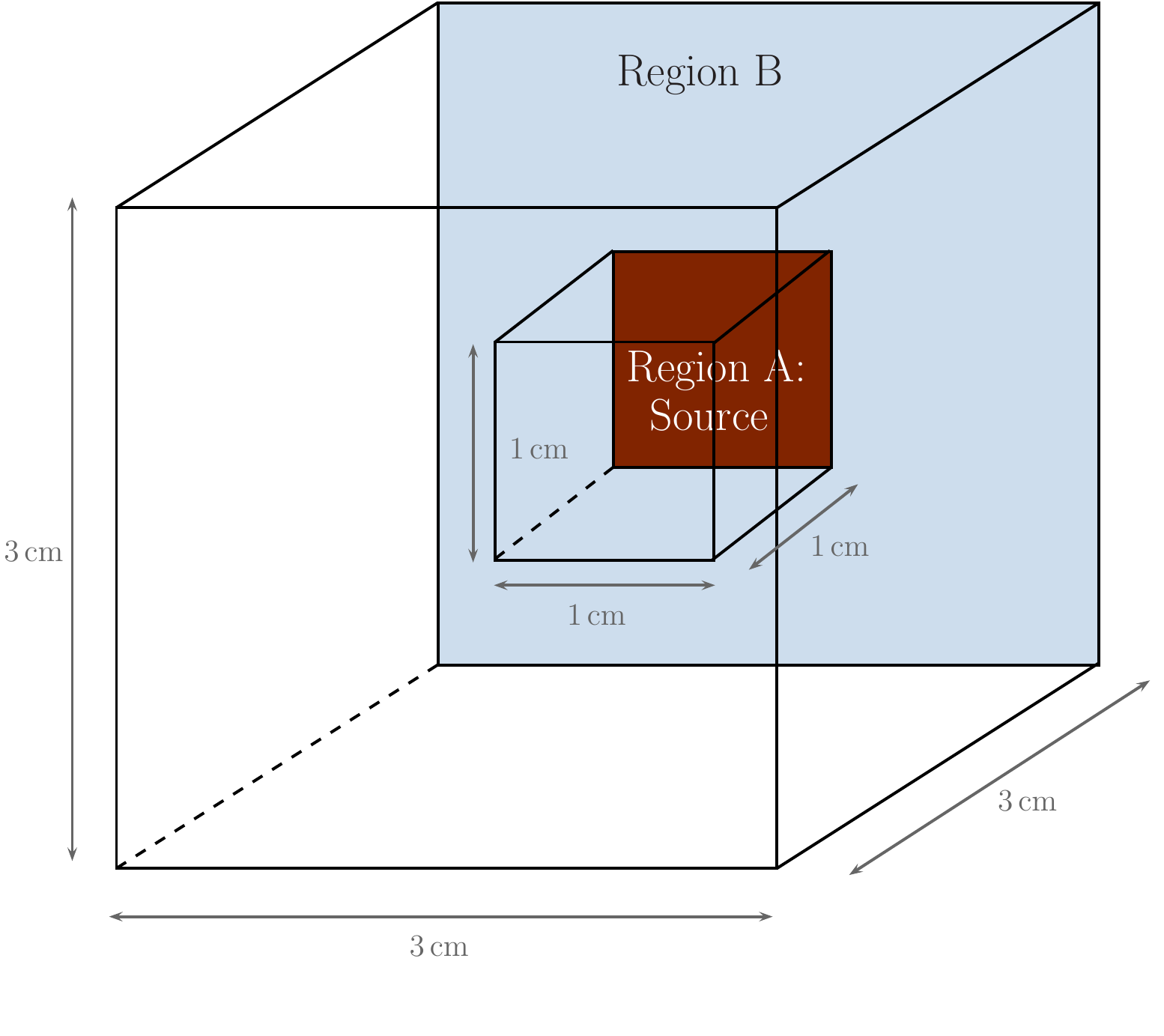}}\\
\caption{Schematic of simple example problems}
\label{fig:eg_probs}
\end{figure}

The two tests problems are based on the same geometry, shown in \fref{fig:eg_probs}, but use different material properties to represent an ``easy'' and ``hard'' problem (see \tref{tab:material_prop} for the material properties). The ``easy'' problem (which we denote as Problem 1 - diffuse) contains a single material and is very diffuse, with a large scatter cross-section. This is a problem in which we would expect multigrid to perform very well in. The ``hard'' problem however (which we denote as Problem 2 - absorbing), contains two different materials, both of which have a relatively low absorbing cross-section. This is away from the diffusive limit of \eref{eq:bte}, while also being multi-material and provides a more difficult problem for our multigrid algorithm to handle. As mentioned above, we use gmsh to mesh the domains in both 2D (resulting in a mesh with 59k elements and 29k nodes) and 3D (resulting in a mesh with 219k elements and 37k nodes). We should note that the coarsenings produced are the same for both Problem 1 and 2.
%~~~~~
\begin{table}[ht]
\centering
\begin{tabular}{ p{2cm} c c c c c c}
\toprule
Problem & \multicolumn{3}{c}{Region A} & \multicolumn{3}{c}{Region B} \\
\midrule
& Source (cm$^{-2}$ s$^{-1}$) & $\Sigma_\textrm{t}$ (cm$^{-1}$) & $\Sigma_\textrm{s}$ (cm$^{-1}$) & Source (cm$^{-2}$ s$^{-1}$) & $\Sigma_\textrm{t}$ (cm$^{-1}$) & $\Sigma_\textrm{s}$ (cm$^{-1}$) \\
\midrule
1 - Diffuse & 1.0 & 10.0 & 10.0 & 0.0 & 10.0 & 10.0 \\
2 - Absorbing & 1.0 & 0.5 & 0.0 & 0.0 & 1.0 & 0.0 \\
\bottomrule  
\end{tabular}
\caption{Material properties for the 2 and 3-dimensional example problems shown in \fref{fig:eg_probs}}
\label{tab:material_prop}
\end{table}
%~~~~~
% ~~~~~~~~~~~~~~~~~
% ~~~~~~~~~~~~~~~~~
\subsubsection{Aggressive coarsening}
To begin, we consider the properties of the resulting mesh after one coarsening. We tested a range of desired agglomerate sizes, in 2D from 4 to 500, while in 3D from 8 to 1600. The results from this are shown in the shaded columns of Tables \ref{tab:2D_size} and \ref{tab:3D_size}, with the average, actual agglomerate size produced by each of the algorithms after one coarsening for the 2D and 3D problems, respectively.
%~~~~~
\begin{table}[ht]
\centering
\begin{tabular}{ p{2cm} a c a c a c}
\toprule
\rowcolor{White}
Desired agglom. size & \multicolumn{2}{c}{\textcolor{Black}{Greedy}} & \multicolumn{2}{c}{\textcolor{foliagegreen}{METIS \cite{karypis_metis-unstructured_1995}}} & \multicolumn{2}{c}{\textcolor{matlabblue}{MGridGen \cite{moulitsas_multilevel_2001}}} \\
\midrule
4 & 3.97 & 1.97 & 4.0 & 1.98 & 3.4 & 2.22\\
8 & 7.01 & 1.59 & 8.0 & 1.57 & 6.4 & 1.75\\
12 & 10.05 & 1.43 & 12.0 & 1.40 & 9.02 & 1.56\\
16 & 13.04 & 1.34 & 16.0 & 1.30 & 12.1 & 1.42\\
20 & 16.05 & 1.28 & 20.0 & 1.25 & 14.9 & 1.35\\
24 & 18.5 & 1.24 & 23.98 & 1.21 & 17.7 & 1.29\\
28 & 21.9 & 1.21 & 27.98 & 1.18 & 20.6 & 1.25\\
100 & 72 & 1.06 & 99.6 & 1.05 & 72.5 & 1.07\\
250 & 182 & 1.03 & 249.9 & 1.02 & 184.5 & 1.03\\
500 & 346 & 1.01 & 497.8 & 1.01 & 356.8 & 1.01\\
\bottomrule  
\end{tabular}
\caption{The effect of changing the desired agglomerate size input to the first coarsening performed by each of the size-based coarsening routines in 2D, for the problem shown in \fref{fig:2D_eg} meshed with 59k triangles. The actual, average agglomerate sizes after the first coarsening are reported in the shaded column. On the lower grids, the desired agglomerate size was fixed as described in \secref{sec:Choosing}. The resulting spatial grid complexities are reported in the unshaded column for each algorithm. The multigrid solve times for the different algorithms and desired agglomerate sizes are shown in \fref{fig:time_vs_complx_2D_both}}
\label{tab:2D_size}
\end{table}
%~~~~~
Across all desired agglomerate sizes, in both 2D and 3D, METIS produces agglomerates that match the desired size the best, with the Greedy algorithm performing the worst. All three algorithms produce varying agglomerate size as the the desired size increases (particularly in 3D). This is even true with METIS, which by design always returns the exact number of partitions requested. For example, in \tref{tab:3D_size}, with a desired agglomerate of 1600 in 3D, the average actual agglomerate size returned by METIS is 1589. We must remember however, that as mentioned in \secref{sec:Cleanup}, a cleanup routine is run after each coarsening, that can modify the agglomerates produced. This is particularly evident when examining the results of the Greedy algorithm. Again the Greedy algorithm returns agglomerates with the exact desired size, but unlike METIS or MGridGen, there is no consideration given to criteria like decreasing surface area/volume ratios or edge-cuts. This forces the cleanup routine to drastically alter the number of agglomerates. In 2D, requesting a desired agglomerate size of 500 from the Greedy algorithm resulted in an average actual size of 346, while in 3D, a desired size of 1600 results in an average size of 620. In general, this does not pose a large problem, as very large agglomerates seem unlikely to produce acceptable grid complexities for our multigrid algorithm (we will discuss this further below), though it highlights the necessity of an effective cleanup routine.
%~~~~~
\begin{table}[ht]
\centering
\begin{tabular}{ p{2cm} a c a c a c}
\toprule
\rowcolor{White}
Desired agglom. size & \multicolumn{2}{c}{\textcolor{Black}{Greedy}} & \multicolumn{2}{c}{\textcolor{foliagegreen}{METIS \cite{karypis_metis-unstructured_1995}}} & \multicolumn{2}{c}{\textcolor{matlabblue}{MGridGen \cite{moulitsas_multilevel_2001}}} \\
\midrule
8 & 7.3 & 2.65 & 8 & 2.48 & 6.9 & 2.48\\
40 & 27.9 & 1.93 & 40 & 1.73 & 32.3 & 1.86 \\
72 & 45.1 & 1.67 & 72 & 1.5 & 59 & 1.60\\
104 & 61.9 & 1.53 & 104 & 1.36 & 85.7 & 1.46\\
136 & 76.1 & 1.44 & 135.0 & 1.29 & 112 & 1.38\\
168 & 91.0 & 1.38 & 168.0 & 1.24 & 138.3 & 1.31\\
200 & 106.5 & 1.33 & 199.9 & 1.21 & 165.7 & 1.26\\
300 & 145.2 & 1.25 & 299.8 & 1.14 & 246.5 & 1.2 \\
600 & 255.5 & 1.14 & 599.7 & 1.08 & 491.19 & 1.1 \\
1200 & 449.4 & 1.06 & 1199.3 & 1.03 & 1007.7 & 1.05 \\
1600 & 620.6 & 1.05 & 1589.4 & 1.02 & 1376.7 & 1.04 \\
\bottomrule  
\end{tabular}
\caption{The effect of changing the desired agglomerate size input to the first coarsening performed by each of the size-based coarsening routines in 3D, for the problem shown in \fref{fig:3D_eg} meshed with 219k triangles. The actual, average agglomerate sizes after the first coarsening are reported in the shaded column. On the lower grids, the desired agglomerate size was fixed as described in \secref{sec:Choosing}. The resulting spatial grid complexities are reported in the unshaded column for each algorithm. The multigrid solve times for the different algorithms and desired agglomerate sizes are shown in \fref{fig:time_vs_complx_3D_both}}
\label{tab:3D_size}
\end{table}
%~~~~~
% ~~~~~~~~~~~~~~~~~
\begin{figure}[p]
\centering
\subfloat[][2D. The node/element ratio on the top grid was $\sim0.5$.]{\label{fig:size_vs_node_ele_ratio_2D_not_reg}\includegraphics[width =0.425\textwidth]{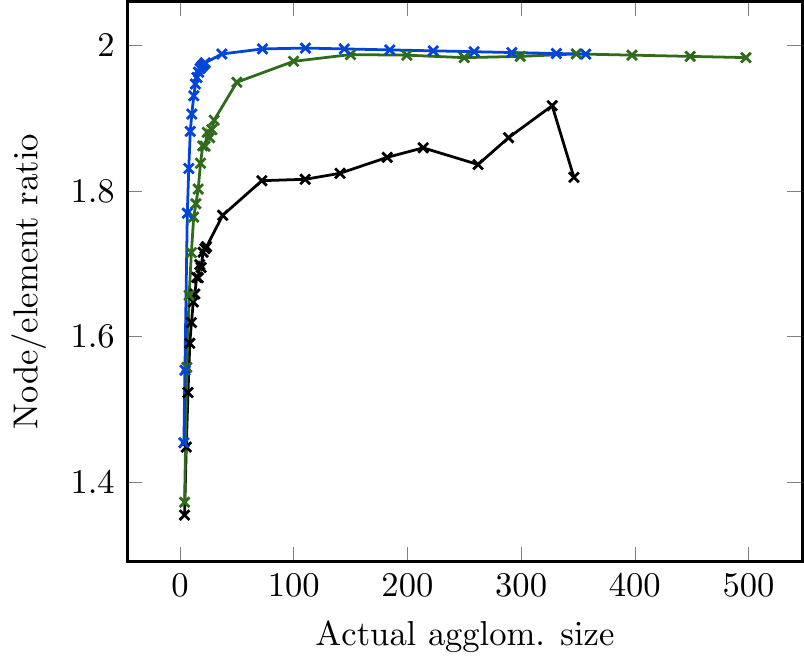}} \hspace{0.5cm}
\subfloat[][3D. The node/element ratio on the top grid was $\sim0.19$.]{\label{fig:size_vs_node_ele_ratio_3D}\includegraphics[width =0.425\textwidth]{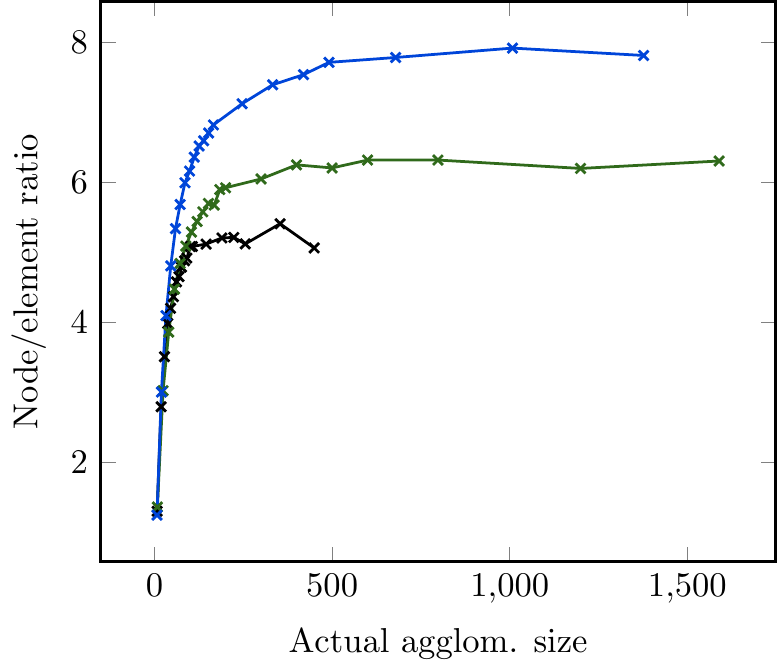}}\\
\caption{Node/element ratio after the first coarsening, given the average actual agglomerate sizes produced after the first coarsening. The black curve denotes the Greedy algorithm, the green {\textcolor{foliagegreen}{METIS \cite{karypis_metis-unstructured_1995}}} and the blue {\textcolor{matlabblue}{MGridGen \cite{moulitsas_multilevel_2001}}}}
\label{fig:size_vs_node_ele_ratio}
\end{figure}
% ~~~~~~~~~~~~~~~~~
% ~~~~~~~~~~~~~~~~~
\begin{figure}[p]
\centering
\subfloat[][2D. The average connectivity on the top grid was $\sim6$.]{\label{fig:size_vs_avg_connect_2D_not_reg}\includegraphics[width =0.425\textwidth]{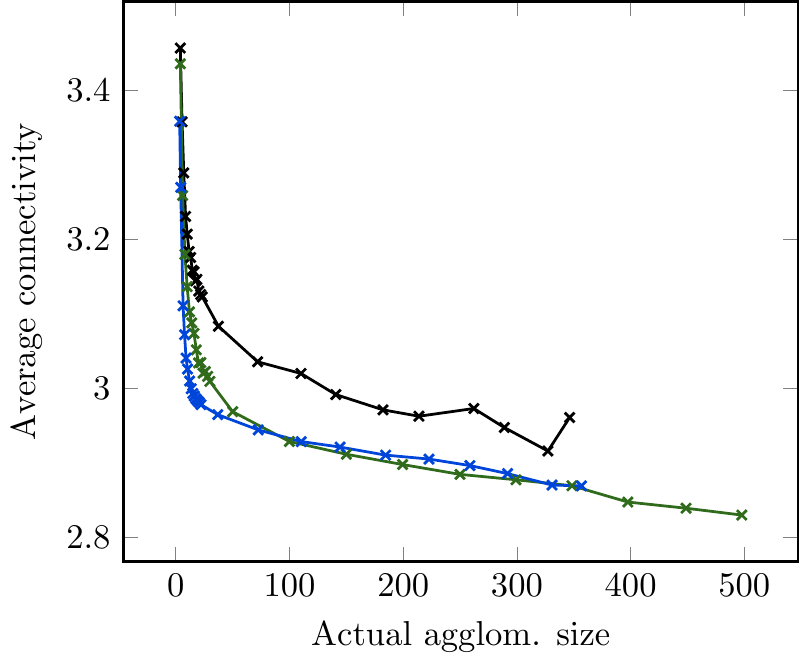}} \hspace{0.5cm}
\subfloat[][3D. The average connectivity on the top grid was $\sim20$.]{\label{fig:size_vs_avg_connect_3D}\includegraphics[width =0.425\textwidth]{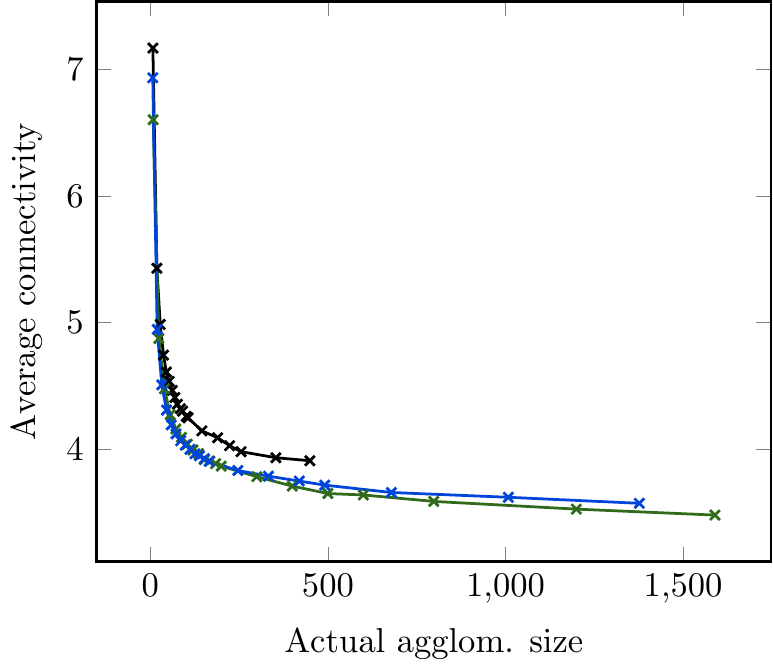}}\\
\caption{Average connectivity after the first coarsening, given the average actual agglomerate sizes produced after the first coarsening. The black curve denotes the Greedy algorithm, the green {\textcolor{foliagegreen}{METIS \cite{karypis_metis-unstructured_1995}}} and the blue {\textcolor{matlabblue}{MGridGen \cite{moulitsas_multilevel_2001}}}}
\label{fig:size_vs_avg_connect}
\end{figure}
% ~~~~~~~~~~~~~~~~~

Regardless of the desired size requested from each of the algorithms or the actual size returned, the connectivity of the mesh after a single coarsening fundamentally changes when compared with the original mesh. This is demonstrated in both Figures \ref{fig:size_vs_node_ele_ratio} and \ref{fig:size_vs_avg_connect}. In 2D, the node/element ratio on the original mesh was approx. $0.5$, meaning there are more elements than nodes (in contrast to a structured quadrilateral mesh). This is reflected in the average connectivity, which is approx. $6$ (larger than on a structured quadrilateral mesh). After a single coarsening however, this is no longer the case. \fref{fig:size_vs_node_ele_ratio_2D_not_reg} shows that the node element ratio becomes greater than $1$, regardless of the actual agglomerate size produced by the coarsening algorithms. Similarly, the average connectivity decreases considerably, as shown in \fref{fig:size_vs_avg_connect_2D_not_reg}. The same is true in 3D, though the difference is even greater. Originally, in 3D, the original mesh had a node/element ratio of approx. $0.19$, though \fref{fig:size_vs_node_ele_ratio_3D} shows this has also becomes greater than $1$ after coarsening, with \fref{fig:size_vs_avg_connect_3D} showing that the average connectivity reduces considerably from approx. $20$ on the original mesh. \fref{fig:size_vs_node_ele_ratio_3D} also shows a significant difference in the node/element ratio between METIS and MGridGen in 3D (from approx. 6 to 8) that is not seen in 2D. Given coarse nodes are defined at the intersection of coarse faces between elements, the optimisation performed \cite{moulitsas_multilevel_2001, Moulitsas2001} by MGridGen of the aspect ratio of coarse elements (the square of the external face area divided by the volume) produces coarse elements closest to spheres and hence affects the node/element ratio. This is more evident in 3D and at large actual agglomerate sizes, as MGridGen can approximate a sphere better with larger numbers of fine elements.

Afer a single coarsening, we see node/element ratios and average connectivities like those on a structured grid. The change in connectivity after the first coarsen suggests that subsequent coarsenings should be performed less aggressively, with the desired agglomerate size smaller than on the top grid. As such, for all size-based agglomeration algorithms, we fix the desired agglomerate size at 4 in 2D and 8 in 3D after the first coarsening. This is a very similar result to that seen in AMG methods as one might expect, where aggressive coarsening are often used on the top multigrid level \cite{trottenberg_multigrid_2000}. Furthermore in 3D, we again reduce the desired agglomerate size to 4 in 3D if there are less than $100$ elements remaining on a coarse grid. 

These simplistic heuristics for the lower grid agglomerate size perform quite well in practice; for example, in the 2D problem using METIS, with a top grid desired agglomerate size varying between 4 and 500, using a desired agglomerate size of 4 on the second grid results in a reduction in number of nodes by a factor of between 3 and 4 respectively. This highlights that we can consider the lower grids in some sense ``structured''. Of course we could try and determine a desired agglomerate size on lower grids automatically, given the node/element ratio and average connectivity on the higher grids, but we do not attempt that in this work. 

%% ~~~~~~~~~~~~~~~~~
%\begin{figure}[hpt]
%\centering
%\subfloat[][2D.]{\label{fig:size_vs_face_nodes_2D_not_reg}\includegraphics[width =0.425\textwidth]{figures/PS/size_vs_face_nodes_2D_not_reg.pdf}} \hspace{0.5cm}
%\subfloat[][3D.]{\label{fig:size_vs_face_nodes_3D}\includegraphics[width =0.425\textwidth]{figures/PS/size_vs_face_nodes_3D.pdf}}\\
%\caption{The ratio of the number of nodes removed from the interior of elements vs the number removed from the faces of elements, after the first coarsening, given the average actual agglomerate sizes. The black curve denotes the Greedy algorithm, the green {\textcolor{foliagegreen}{METIS \cite{karypis_metis-unstructured_1995}}} and the blue {\textcolor{matlabblue}{MGridGen \cite{moulitsas_multilevel_2001}}}. The grey dotted line denotes a ratio of 1.}
%\label{fig:size_vs_face_nodes}
%\end{figure}
%% ~~~~~~~~~~~~~~~~~

\subsubsection{Ideal agglomerate size on top grid}
\label{sec:Resulting}
Given that we have fixed the desired agglomerate size on the lower grids, we can now examine the effect of changing the desired agglomerate size on the top grid only. The unshaded columns of Tables \ref{tab:2D_size} and \ref{tab:3D_size} show the spatial grid complexity in 2D and 3D, respectively, as the desired agglomerate size is changed on the top grid (as mentioned previously, the shaded column shows the actual agglomerate sizes); lower grids have the fixed desired size as discussed above. We can see that the spatial grid complexities are very large for small agglomerate sizes, and decrease as the desired agglomerate size is increased. This is to be expected, as we can remove more fine nodes from the interiors of large elements. The interesting result from Tables \ref{tab:2D_size} and \ref{tab:3D_size} is how large agglomerates must be made (approx. $20$ in 2D, $300$ in 3D) in order to recover grid complexities that approach the ``ideal'' complexity on a structured quadrilateral grid (4/3 in 2D and 8/7 in 3D), regardless of the coarsening algorithm chosen. This is again due to the increased connectivity on an unstructured mesh.

\begin{figure}[p]
\centering
\subfloat[][Problem 1 - diffuse.]{\label{fig:time_vs_complx_2D_diff}\includegraphics[width =0.425\textwidth]{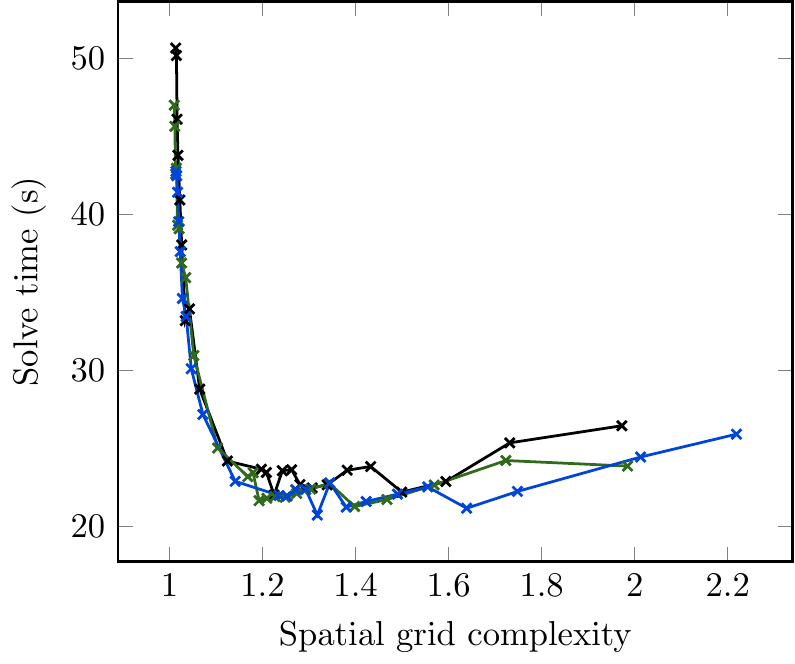}} \hspace{0.5cm}
\subfloat[][Problem 2 - absorbing.]{\label{fig:time_vs_complx_2D}\includegraphics[width =0.425\textwidth]{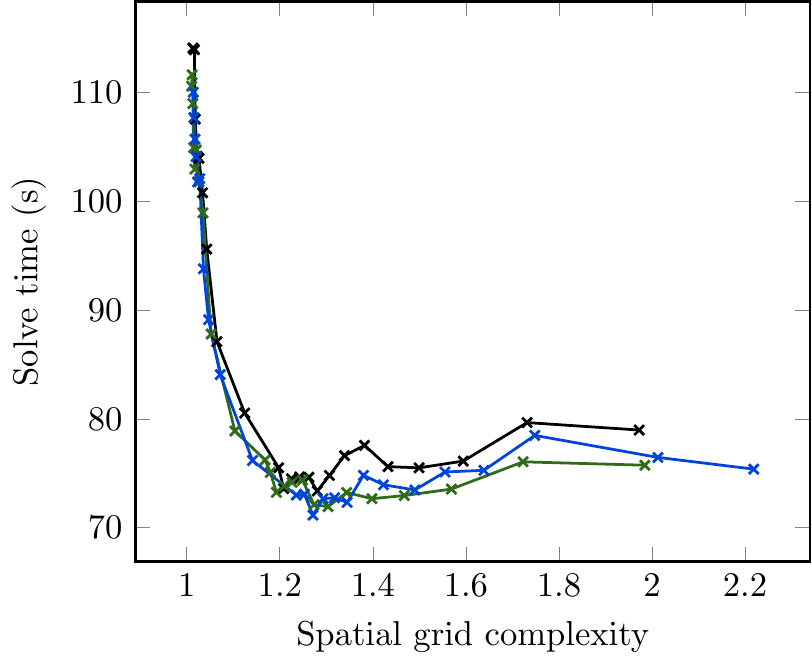}}\\
\caption{Solve time given different spatial grid complexities for the size-based coarsening algorithms in 2D, for the problem shown in \fref{fig:2D_eg} meshed with 59k triangles. The different spatial grid complexities were produced by changing the desired agglomerate size for the first coarsening only (some of these desired agglomerate sizes and resulting spatial grid complexities are given in \tref{tab:2D_size}). On the lower grids, the desired agglomerate size was set to 4. The black curve denotes the Greedy algorithm, the green {\textcolor{foliagegreen}{METIS \cite{karypis_metis-unstructured_1995}}} and the blue {\textcolor{matlabblue}{MGridGen \cite{moulitsas_multilevel_2001}}}}
\label{fig:time_vs_complx_2D_both}
\end{figure}

\begin{figure}[p]
\centering
\subfloat[][Problem 1 - diffuse.]{\label{fig:time_vs_complx_3D_diff}\includegraphics[width =0.425\textwidth]{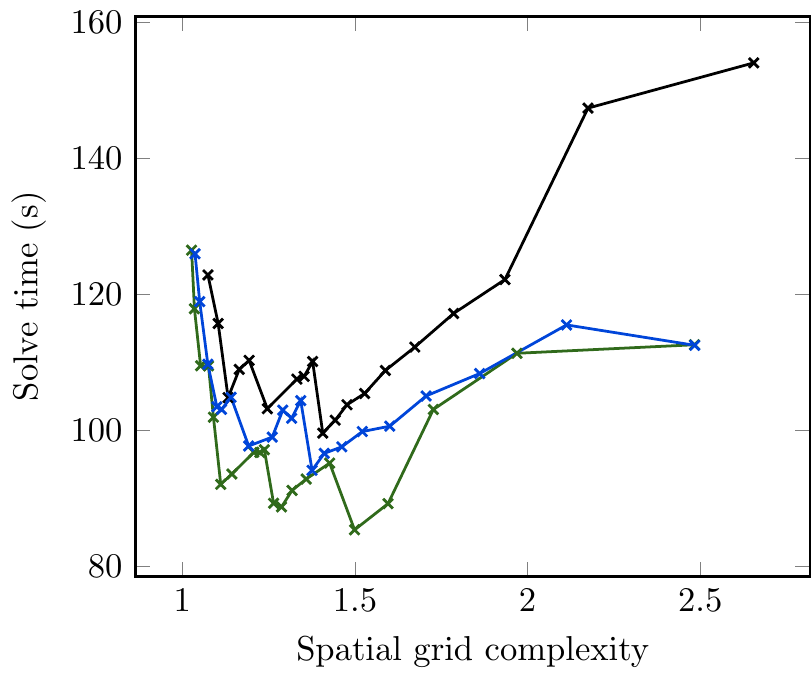}} \hspace{0.5cm}
\subfloat[][Problem 2 - absorbing.]{\label{fig:time_vs_complx_3D}\includegraphics[width =0.425\textwidth]{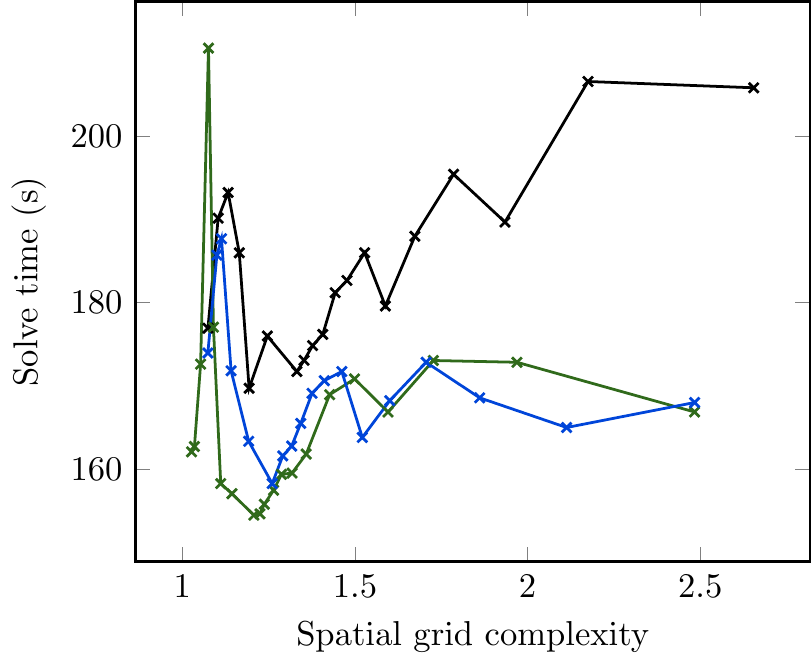}}\\
\caption{Solve time given different spatial grid complexities for the size-based coarsening algorithms in 3D, for the problem shown in \fref{fig:3D_eg} meshed with 223k tetrahedrons. The different spatial grid complexities were produced by changing the desired agglomerate size for the first coarsening only (some of these desired agglomerate sizes and resulting spatial grid complexities are given in \tref{tab:3D_size}). On the lower grids, the desired agglomerate size was set to 8. The black curve denotes the Greedy algorithm, the green {\textcolor{foliagegreen}{METIS \cite{karypis_metis-unstructured_1995}}} and the blue {\textcolor{matlabblue}{MGridGen \cite{moulitsas_multilevel_2001}}}}
\label{fig:time_vs_complx_3D_both}
\end{figure}

In a general sense, we also find the smallest solve time using our multigrid algorithm when we target close to ``ideal'' grid complexities from structured grids. Figures \ref{fig:time_vs_complx_2D_both} and \ref{fig:time_vs_complx_3D_both} show the solve times of our iterative method, as the desired agglomerate size is increased (plotted against the spatial grid complexity from Tables \ref{tab:2D_size} and \ref{tab:3D_size}) for both Problem 1 and 2, as defined in \secref{sec:Choosing}.

\fref{fig:time_vs_complx_2D_diff} shows that in 2D, for the diffusive problem, the solve time for all coarsening algorithms is almost constant for grid complexities greater than $1.2$ (corresponding to an agglomerate size of less than $40$). As long as the agglomerate size does not become too large or too small, the actual desired agglomerate size on the top grid does not matter. For the absorbing problem in 2D, shown in \fref{fig:time_vs_complx_2D}, this is also true with a minimum for all coarsening algorithms occurring at a grid complexity of approx. $1.25$, corresponding to a desired agglomerate size of around 24.

In 3D, \fref{fig:time_vs_complx_3D_diff} shows that for the diffusive problem, we do not see the constant solve time across a range of agglomerate sizes as in \fref{fig:time_vs_complx_2D_diff}. A minimum occurs at a grid complexity of between $1.35-1.5$, which corresponds to a desired agglomerate size of around 72-200. We can see a number of large jumps in the solve times for each coarsening algorithm; these correspond to an increase in the iteration count. For example with a grid complexity of 1.5, METIS takes 8 iterations, at 1.42 it takes 9, and the transition to 10 iterations occurs at the grid complexity of 1.24.

For the absorbing problem shown in \fref{fig:time_vs_complx_3D}, the solve time is higher than the diffusive as would be expected, the minimum solve time occurs at a grid complexity of around 1.25. We do not see such distinct ``jumps'' in solve time shared across all the coarsening algorithms when compared to the diffuse case in \fref{fig:time_vs_complx_3D_diff}, as the iteration count is larger (and hence the impact of performing a single extra iteration is smaller) and we see a more constant iteration count across a range of grid complexities, e.g., using METIS takes 16 iterations between a grid complexity of 1.2--1.5. 

In general, over all desired agglomerate sizes for both problems, we can see from Figure \ref{fig:time_vs_complx_3D_both} that the Greedy algorithm tends to perform significantly poorer in 3D than either METIS or MGridGen. We can also see that METIS consistently result in the smallest solve times across almost all grid complexities. The ideal grid complexity in both these problems in 3D is therefore between around 1.25--1.5; this is larger than might be expected when comparing to structured geometric coarsenings. As mentioned previously, larger agglomerates remove increasingly larger fractions of nodes from their interior, decreasing the locality/accuracy of our simple interpolation operator. 

The results from the 3D cases show that for all of the three coarsening algorithms, the solve time can vary considerably, even in regions with very similar grid complexities. The most obvious of these variations occurs in the absorbing problem, when using METIS with a target agglomerate size on the top grid of 600, giving a grid complexity of 1.08 (see \fref{fig:time_vs_complx_3D}). The iteration count increases to 23 and the solve time is over 200 seconds. This is due to a small number of poorly shaped coarse agglomerates, which are much more likely in 3D, particularly as the target agglomerate size grows (and hence the grid complexity decreases) and are almost impossible to completely eliminate. Thankfully this doesn't occurs when using METIS or MGridGen when in a region of ``sensible'' grid complexities between 1.2--1.5, as we see fairly smooth changes in the solve time.

Given these results, the optimal spatial grid complexity for these problems seems to be at least close to that given by a coarsening on structured grids. This is perhaps an obvious result, but what should be emphasised is that this is fundamentally different to targeting the agglomerate size returned by coarsening a structured grid. Indeed, the agglomerate size necessary to generate acceptable grid complexities is very large, around 24 in 2D and greater than $100$ in 3D. This is at odds with how most element-based coarsening algorithms are designed (e.g., the \citet{jones_amge_2001} and \citet{kraus_agglomeration-based_2004} algorithms), as when used in an AMGe setting (as opposed to our purely geometric multigrid), assumptions are often made about the connectivity of all fine nodes in an agglomerate with coarse nodes, preventing the use of large agglomerates. This leads to the larger grid complexities typical of an AMG approach. This is seen in \citet{wabro_amgecoarsening_2006}, who use the \citet{jones_amge_2001} and RGB algorithms, along with the METIS and MGridGen libraries to perform coarsenings. Even when using the size-based algorithms METIS and MGridGen, they specify a very small desired agglomerate size (4 or 5) in both 2D and 3D, resulting in large grid complexities (typically between 1.8 and 5.3 in 3D).

Whether larger agglomerates could be used in an AMGe setting to reduce the grid complexities, while retaining acceptable convergence is an open question, though for our purely geometric approach small agglomerates are not acceptable. As such, for all the size-based coarsening algorithms shown below, we fix the desired agglomerate size at 24 on the top grid and 4 on lower grids in 2D; in 3D we set 168 on the top grid, with 8 on the lower grids unless there is fewer than 100 elements, in which case we further reduce the aggressivity and target 4. This gives close to the minimum runtime we found for each case with each of the different size-based methods and makes comparisons in the next section simpler. 

\subsection{Comparison between algorithms}
\label{sec:performance}
%~~~~~
\begin{table}[ht]
\centering
\begin{tabular}{ c p{3cm} c c c c c c c}
\toprule
\rowcolor{White}
& & \textcolor{darkpurple}{Jones \cite{jones_amge_2001}} & \textcolor{darksalmon}{Kraus \cite{kraus_agglomeration-based_2004}} & \textcolor{fireenginered}{RGB \cite{wabro_amgecoarsening_2006}} & \textcolor{deludedorange}{Node} & \textcolor{Black}{Greedy} & \textcolor{foliagegreen}{METIS \cite{karypis_metis-unstructured_1995}} & \textcolor{matlabblue}{MGridGen \cite{moulitsas_multilevel_2001}} \\
\midrule
\multirow{1}{*}{2D}& Grid complexity & 1.78 & 1.80 & 1.61 & 1.77 & 1.24 & 1.21 & 1.29 \\
& Multigrid setup time & 6.9 & 7 & 4.7 & 5.5 & 3.6 & 4.3 & 3.9 \\
\rowcolor{light}
& Prob. 1 runtime & 34.6 & 35.2 & 34.8 & 32.9 & 31.1 & 30.04 & 30.3 \\
& No. iterations & 11 & 11 & 13 & 11 & 14 & 13 & 13 \\
\rowcolor{light}
& Prob. 2 runtime & 90.6 & 89.7 & 91.5 & 85.7 & 82.2 & 82.1 & 80.6 \\
& No. iterations & 38 & 37 & 42 & 37 & 44 & 44 & 42 \\[0.5cm]
\multirow{1}{*}{3D}& Grid complexity & 3.22 & 2.54 & 2.29 & 2.41 & 1.38 & 1.24 & 1.32 \\
& Multigrid setup time & 106 & 95.8  & 60 & 63.5 & 43.8 & 31.1 & 34.8 \\
\rowcolor{light}
& Prob. 1 runtime & 376.9 & 296.6 & 243.6 & 261.5 & 230.9 & 176.2 & 195.7 \\
& No. iterations & 8 & 8 & 9 & 8 & 10 & 10 & 10 \\
\rowcolor{light}
& Prob. 2 runtime & 432 & 370.4 & 283.3 & 335.1 & 296.1 & 236.5 & 255.6 \\
& No. iterations & 11 & 13 & 12 & 13 & 16 & 16 & 16 \\
\bottomrule  
\end{tabular}
\caption{The grid complexities and runtimes (s) for the two different problem types described in \tref{tab:material_prop} and \fref{fig:eg_probs}, in 2D and 3D, given different coarsening algorithms. A description of how the desired agglomerates sizes were chosen for the Greedy, METIS and MGridGen algorithms is given in \secref{sec:Choosing}.}
\label{tab:compare_size}
\end{table}
%~~~~~
Now that we have fixed the desired agglomerate sizes on each level for the size-based algorithms, \tref{tab:compare_size} compares the total runtime of all seven coarsening algorithms, in 2D and 3D for the two example problems described in \secref{sec:Choosing}. In 2D, there is no more than a 20\% difference in runtime between the different algorithms, in both the diffusive and absorbing problem. The multigrid setup time is also fairly consistent, however both the Kraus and Jones algorithms requiring up to twice the time to coarsen. 

The results in 3D show far more difference between the algorithms, with METIS and MGridGen performing the best, both in terms of total runtime and multigrid setup time, with Jones performing the worst. The poor performance of Jones is to be expected, given the Kraus algorithm is a modified version of the Jones algorithm designed to improve the coarsenings in 3D, and we can see that Kraus does perform better. We should also note that METIS/MGridGen required much less ``cleanup'' (see \secref{sec:Cleanup}) than the other algorithms; this becomes very noticeable in 3D, particularly for the simple Node and Greedy algorithms.

For the diffusive problem in 3D, coarsening with METIS gives a total runtime over half that when compared with the Jones algorithm. The multigrid setup time for Jones is also three times that of METIS. The Jones and Kraus algorithms were both designed to cost $\mathcal{O}(n)$ in the number of faces (or edges for Kraus), however these algorithms were constructed to try and replicate the same element ratio as on structured grids for AMGe applications. For example, the Jones algorithm results in an average of 4.91 fine elements in each agglomerate after the first coarsening, with Kraus giving 9.5; this is in comparison to the 168 we determined to be optimal for the size-based algorithms in \secref{sec:Resulting}. 

The multigrid setup time shown in \tref{tab:compare_size} clearly depends on the efficiency of our implementation of each coarsening algorithm (especially given our use of a well-used library like METIS), coarse topology creation, cleanup step, etc., but small agglomerate size results in many coarse faces and profiling reveals that a significant component of the setup time is dependent on the number of faces. This is both when when performing the coarsening (given the $\mathcal{O}(n)$ behaviour) and also when creating the coarse topology; the original fine mesh has 782k faces, with the Jones algorithm resulting in 293k coarse faces and METIS giving only 16k. This helps explain the significant difference in setup time and this trend would remain even if our implementations were improved. 

\tref{tab:compare_size} shows a significant difference in the grid complexities between the different algorithms in 3D, with 3.22 for Jones, and 1.24 for METIS. The larger grid complexities lead to fewer iterations (in the absorbing problems Jones took 11 iterations, with METIS at 16), but as is common in multigrid methods the extra expense of smoothing and computing matvecs on the larger coarse grids does not result in a smaller runtime with such large grid complexities. 

We could examine the operator complexity directly to further quantify the cost of computing matvecs on our lower grids (though it is difficult to compute in our matrix-free setting), but given our simple interpolation (and equivalent FEM rediscretisation on coarse grids), the cost is linearly related to the number of coarse elements on a grid and quadratic in the number of nodes on a coarse element (in the average). METIS benefits from having approx. $34$ times fewer coarse elements on the second grid when compared to Jones, but the smaller coarse elements in Jones give a node/element ratio of 0.88, in comparison to METIS with 5.37, giving 6.3 times more nodes per element with METIS. This balance results in matvecs on the second grid with Jones measured at 5.4 times the cost of that with METIS. 

This shows that in our geometric setting it is worth agglomerating together many elements to reduce the grid complexity, even though this comes at the increased (quadratic) cost of operating on more nodes on each coarse element, given that fine nodes are eliminated from the interior of large coarse elements quicker than on coarse faces (see \secref{sec:Coarse node selection}). These trends also hold for the other coarsening algorithms and these combination of factors help explain the strong relationship between grid complexity and runtime in 3D seen in \tref{tab:compare_size}.

Of course we should note that the Jones, Kraus, RGB and Node algorithms are behaving as expected and would likely perform very well in AMGe settings when compared to using Greedy/METIS/MGridGen with such large target agglomerate sizes. Indeed comparing the total runtimes for METIS for example, with a target agglomerate size of 5 (in an attempt to match the result from the Jones algorithm, which results in an average of 4.91 fine elements per agglomerate) in the 3D absorbing problem results in a grid complexity of 2.57, a multigrid setup time of 62 seconds and a total runtime of 320 seconds, which is much closer to that of the Jones algorithm. 
% ~~~~~~~~~~~~~~~~~
% ~~~~~~~~~~~~~~~~~
% ~~~~~~~~~~~~~~~~~
% ~~~~~~~~~~~~~~~~~
%~~~~~~~~~~~~~~~~~~~~~~~~~~~~~~~~~~~~
\section{Conclusions}
This paper has examined the performance of seven different element coarsening algorithms when used as part of a geometric multigrid on unstructured grids for radiation transport applications. Unstructured tri/tet meshes naturally feature nodes that are connected to more elements than their structured counterparts, regardless of the quality of the mesh and hence using coarsening algorithms that target agglomerate sizes that mimic structured coarsenings instead of targetting grid complexities leads to an inefficient geometric multigrid method, as would be expected.

We found that coarsening algorithms which allowed the easy specification of target agglomerate sizes, namely using METIS or MGridGen the most suited to unstructured geometric multigrid, as we could easily force the coarsening on the top grid to be aggressive by agglomerating many fine elements, while decreasing the aggressiveness on lower grids. The smallest runtimes resulted from grid complexities that were close to those found in structured multigrids in 2D (around 1.25), while slightly larger in 3D (between 1.25--1.5).

We showed that ``simple'' coarsening strategies like agglomerating all elements connected to a node, or a naive size-based search can result in acceptable coarse grids with moderate runtimes (particularly in 2D), though we found that these algorithms result in poor agglomerates that required considerable ``cleanup'' to form sensible coarse grids. These methods can be very useful however in specialised applications, like GPU-based coarsenings where improving the speed/parallelism of the coarsening with a simple algorithm may result in smaller total runtimes even if the resulting multigrid is less than optimal. 

Coarsening algorithms used in AMGe-type methods, like the Jones, Kraus or RGB algorithms performed the worst in this application, as their desire to limit the coarse element size in order to provide AMGe-type interpolation meant that in a purely geometric setting with simple operators the resulting agglomerates were too small, giving large grid complexities and hence runtimes. 

The key message from this work is simple and matches common knowledge from AMG methods: if an element agglomeration algorithm is used on an unstructured mesh to form coarse grids for a geometric multigrid, use an aggressive coarsening on the top grid (agglomerate 24 elements in 2D, 170--200 in 3D). This results in a second grid with fundamentally different connectivity and the subsequent coarsenings can be less aggressive (4 in 2D and 8 in 3D). Furthermore given the prevalence in existing scientific software, using METIS to produce this coarsening is thankfully an excellent choice and almost always resulted in the smallest runtimes and multigrid setup times while producing ``nice'' coarsenings that required little cleanup. The extra optimisation performed by MGridGen to produce agglomerates shaped close to circles/spheres were visible (particularly in 3D) but did not seem to result in smaller runtimes when compared to METIS. 

The recommendation of using METIS to perform coarsenings has the added benefit of making migration to an AMGe-type method very simple, as the targetted agglomeration size can be easily changed. If new element coarsening algorithms are to be developed, this work suggests that allowing the size of agglomerates to be easily changed is a key feature that would also help extend the applicability beyond AMGe-type multigrids. 
%~~~~~~~~~~~~~~~~~~~~~~~~~~~~~~~~~~~~
\section*{Acknowledgements}
The authors would also like to acknowledge the support of the EPSRC through the funding of the EPSRC grants EP/P013198/1, EP/R029423/1 and EP/M022684/2.
%~~~~~~~~~~~~~~~~~~~~~~~~~~~~~~~~~~~~
%~~~~~~~~~~~~~~~~~~~~~~~~~~~~~~~~~~~~
%~~~~~~~~~~~~~~~~~~~~~~~~~~~~~~~~~~~~



%% References

\section*{References}
\bibliographystyle{model1-num-names}
\bibliography{bib_library}






% \end{thebibliography}

\end{document}